\input amstex
\documentstyle{amsppt}

%Begin Macros
\define \cal {\Cal} 
\define \p {\cal P}
\define \Q {{\Bbb Q}}
\define \Z {{\Bbb Z }}
\define \R {{\Bbb R }}
\define \HH {{\Bbb H }}
\define \C {{\Bbb C} }
\define \bC {\C }
\define \PP {{\Bbb P }}
\define \N {{\Bbb N }}
\define \G {\cal G}

\define \La {\Lambda}
\define \LC {\Lambda_\C}
\define \LCs {\Lambda_{\C^*}}
\define \LR {\Lambda_\R}
\define \LQ {\Lambda_\Q}

\define \pf {\noindent {\bf Proof. }}

\define \WU {  W_{U(1)}  }
\define \Ci {C^\infty}
\define \del {\nabla}
\define \cL  {{\cal L}}
\define \cV  {{\cal V}}

\define \Om #1 {\Omega_X^{#1}}
\define \OmD #1 {\Omega_{\bar X}^{ #1 }(\log D)}
\define \OmL  {\Omega_{\bar X}^\bullet (\log D) \otimes {\bar L} }
\define \OmDL #1  {\Omega_{\bar X}^ { #1 } (\log D) \otimes {\bar L} }
\define \OmV  {\Omega_{\bar X}^\bullet (\log D) \otimes {\bar V} }
\define \OmDV #1  {\Omega_{\bar X}^ { #1 } (\log D) \otimes {\bar V} }
\define \Omdt {{\Omega_X^\bullet (log D)}}
\define \dt {\bullet}
\define \lblarr #1  {{\mathop {\longrightarrow} \limits^{#1}}}
\define  \A #1 {A_{\bar X}^{#1} (\log D) \otimes {\bar L}} 
\define \AV #1 {A_{\bar X}^{#1} (\log D) \otimes {\bar V}}
\define \LOm #1 {A_{X, (2)}^ {#1} (V) }

\define \calA {{\cal A}}
\define \calB {{\cal B}}

\define \> {\quad}
% End Macros

\refstyle{A}
\NoBlackBoxes

%%%%%%%%%%%%%%
\topmatter
\title Geometry of Cohomology Support Loci for Local Systems 
I.\endtitle
\author Donu Arapura \endauthor
 \address Department of Mathematics, Purdue University,
  West Lafayette, Indiana 47907 \endaddress
\email dvb\@math.purdue.edu\endemail
\thanks Author partially supported by NSF. \endthanks

\toc
\widestnumber\head{III}
\head Introduction \endhead
\head I. Preliminaries on local systems\endhead
 \subhead 1. Connections \endsubhead
  \subhead 2. Cohomology support loci \endsubhead 
\head II. Exponentials of Hodge Structures\endhead
  \subhead 1. 1-Hodge structures \endsubhead  
  \subhead 2. The subtorus theorem \endsubhead  
\head III.  $\C$-Mixed Hodge theory\endhead 
 \subhead 1. $\C$-mixed Hodge structures\endsubhead 
  \subhead 2. $\C$-Hodge complexes\endsubhead 
 \subhead 3. Key constructions \endsubhead  
\head IV. Hodge Decomposition for Higgs Bundles\endhead 
  \subhead 1. Real analytic log complex\endsubhead 
  \subhead 2. Main theorem\endsubhead 
\head V.  Structure of Rank One Loci\endhead 
\head {} References \endhead
\endtoc

\endtopmatter
%%%%%%%%%%%%%%%%%%%%%%%
\document

\head Introduction.\endhead Let $X$ be a connected topological space which is
homotopic to a finite simplicial complex. The fundamental group of $X$
can be studied via  its representations, or equivalently in
terms of locally constant sheaves on $X$. Consider the set
 $H^1(\pi_1(X),GL_n(\C))_{ss}$ of conjugacy classes of all semisimple 
 $n$ dimensional representations. This set has the structure of a (possibly
reducible) complex affine algebraic variety.  This variety contains
some canonically defined subsets (in fact subvarieties): the cohomology
support loci  defined by
$$\Sigma_m^k(X, GL_n(\C)) 
= \{ \rho \in H^1(\pi_1(X), GL_n(\C))_{ss} 
| \, dim \, H^k(X, (\C^n)_{\rho})  \ge m \}, $$
where $(\C^n)_{\rho}$ is the locally constant sheaf with monodromy
representation $\rho$. These sets  
depend only on the homotopy type of $X$. 

  The primary goal of this work
is to obtain a better geometric understanding of these sets when $X$
is the Zariski open subset of a compact K\"ahler manifold (such as
a smooth quasiprojective variety); these results lead to new
restrictions on the homotopy type of such spaces. 
In this paper, we will primarily
study the rank one case, and leave the higher rank case for the sequel.

 Let $\bar X$ be a compact K\"ahler manifold, $D \subset \bar X$ a divisor
with normal crossings, and $X= \bar X - D$. In chapter V, we obtain
the basic structure theorem of
$\Sigma^k(X)= \Sigma_1^k(X,\C^*)$.

\proclaim {Theorem} If either $X = \bar X$ or $H^1(\bar X, \C) = 0$ then
for each integer $k \ge 0$, there exists a finite number of unitary characters
$\rho_i \in H^1(X, \C^*)$, and holomorphic maps into complex  tori (see below)
$f_i\colon X \to T_i$ such that 
$$\Sigma^k(X) = \bigcup_i \rho_if_i^*H^1(T_i, \C^*).$$
\endproclaim

When $X$ is compact this follows from the work of Green and Lazarsfeld 
\cite{GL2}. In \cite{A},  we outlined an alternative argument using Higgs bundles 
techniques and later Simpson \cite{S2} gave somewhat sharper results along these
lines for  $X$ a projective variety. 
Unfortunately we have not yet overcome all the technical
hurdles in the noncompact case,
  so we have the theorem only under the above restrictive hypothesis.
However, using slightly different arguments, we obtain the following 
generalization of a result of Beauville \cite{B1,B2} without any such restrictions:

\proclaim {Theorem} There exists a finite number of torsion characters 
 $\rho_i \in H^1(X, \C^*)$,
unitary characters $\rho_j'$
and surjective holomorphic maps onto smooth curves $f_i\colon X\to C_i$ such
that
$$\Sigma^1(X) = 
\bigcup_i \rho_i f_i ^*H^1(C_i, \C^*) \,\cup \, \bigcup_j \{\rho_j'\}$$
\endproclaim

The $T_i$ of the first theorem would correspond to the generalized
Albanese varieties of the $C_i$.
The second theorem yields a number of interesting corollaries. For example,
it will be shown that the fundamental group of $X$ determines the
number of surjective holomorphic maps, with connected fibers, of $X$
onto curves with negative Euler characteristic.

The first theorem is deduced as a corollary of another theorem which
shows that $\Sigma^k(X)$
is of exponential Hodge type. This  will mean that there are a finite number of
sub-mixed Hodge structures $\La_i \subset H^1(X,\Z)$ and elements $r_i\in \sqrt{-1}H^1(X,\R)$
such that under the exponential map $\exp\colon H^1(X,\C)\to H^1(X,\C^*)$ we have:
$$\Sigma^k(X) = \bigcup_i \exp(\La_i \otimes \C + r_i).$$

In chapter II, a criterion is given for a subset 
$\Sigma \subseteq H^1(X,\C^*)$  to be of  exponential Hodge type. When 
$X$ is compact this boils down to 2 conditions: the set should
be Zariski closed, and if $\rho\exp(\theta ) \in \Sigma$ with $\rho 
\in H^1(X,U(1))$ and $\theta \in H^0(X,\Omega_X^1)$ then $\rho\exp( t\theta) \in \Sigma $ 
for all real $t$ close to 1. The Zariski closedness of $\Sigma^k(X)$
(and more generally $\Sigma_m^k(X,Gl_n(\C))$)
is relatively straight forward and is dealt with chapter I. 
The second condition can be checked using 
the results established in chapters II, 
III and IV (explained below) which imply that cohomology is invariant under 
rescaling $\theta$.

Let  $(V,\del)$ be a pair  consisting of a holomorphic vector
bundle on $X$ with a flat unitary connection. It
extends to a pair  $(\bar V,\bar \del)$ consisting of a vector bundle on $\bar X$ and a
singular connection in a natural fashion.
 A (log) Higgs structure on $\bar V$ is a holomorphic map
 $\theta\colon \bar V \to \Omega_X^1(log D)\otimes \bar V$ satisfying 
$\theta\wedge\theta=0$. 
 Given a Higgs structure, the graded sheaf
 $\OmV$ can be made into a complex using
 either $\bar \del+\theta$ or just $\theta$ as the differential.
 The operator $\bar \del + \theta$ defines a new flat (nonunitary) 
 connection on $V$, so we obtain a locally constant sheaf of flat 
 sections $\cV_{\theta}$ on $X$.
 If 
 $\theta \in H^0(\Omega_{\bar X}^1\otimes End(\bar V))$ is a Higgs field without poles, then 
we establish isomorphisms of hypercohomology groups:
$${\Bbb H}^\dt(X,\cV_{\theta})\cong{\Bbb H}^\dt(\OmV ;
\bar \del+\theta) 
\eqno{(1)}$$
$$ {\Bbb H}^\dt(\OmV ; \bar \del+\theta) \cong 
{\Bbb H}^\dt(\OmV ; \theta) \eqno{(2)}.$$ 
The last result, which is theorem 2.1 of chapter IV,
 is closely related to a theorem of Simpson \cite{S1, 2.2}. 
Simpson's theorem implies that for compact $X$, an analogue of (2) holds where
the left side would be replaced by 
$${\Bbb H}^\dt(\OmV ; \bar \del+\theta+\bar\theta).$$
As a  corollary of these results, we find that the cohomology of 
$\cV_{\theta}$ is invariant under the dilation $\theta\mapsto t\theta$. 
And this is a crucial step in the proof of the above theorems as we have
already indicated. It would be reasonable to conjecture that the isomorphism (2)
holds even when the Higgs field $\theta$ has poles. As evidence for 
this, we prove  in theorem 2.4 of chapter IV that 
(2) also holds when $\bar V$ is a line bundle 
and $\theta$ is of type $(1,1)$.

 In the simplest case where $X$ is compact and $V=L$ is a line bundle,
  $(2)$  can be proved by showing that 
the spectral sequences
  associated to the stupid filtrations abutting to either side
 degenerate at
$E_2$ and the $E_2$ terms are easily see to coincide.
We know of two proofs of these degeneration results. The first proof
which uses the $\partial\bar\partial$-lemma (compare \cite{GL1, 3.7}) is quite short, 
however it does not  generalize
 to the ``log" setting in any obvious way. The second approach,
which is pursued here, is
modeled on ideas from  mixed Hodge theory. The two spectral
sequences
are realized as direct summands of a bigger  spectral sequence
associated to a $\C$-mixed Hodge complex, and the arguments of
Deligne \cite{D} can be adapted to prove degeneration. The
formalism is worked out in chapter III.

 I would like to express my thanks to H\'el\`ene Esnault and
Eckart Viewheg for their hospitality at Essen where much of the
work for the third and fourth chapters was done.

\subhead Notation and terminology\endsubhead

 If $\La$ and $A$ are abelian groups we write $\La_A = \La \otimes_\Z A$.

Unfortunately
the word ``torus" will be used in many different ways in this paper, so  let us set forth some 
conventions  to distinguish among them.
A {\it real torus} is a real Lie group isomorphic to a product of circles. A {\it compact complex
torus} is a complex Lie group isomorphic, as  a real Lie group, to a product of circles. An
{\it affine torus} is a complex algebraic group isomorphic to a product of $\C^*$'s. A {\it
complex torus} is complex Lie group which is an extension of a compact complex torus by an affine
one. Using the exponential map, one can see that any complex torus is
 a quotient of a complex vector space by a discrete subgroup which
contains a basis of the vector space. 

Manifolds will be assumed to be connected unless stated otherwise. By  a 
curve, we will mean a compact connected complex curve with finitely many
points removed. $\Omega_{\bar X}^p$ will denote the sheaf of holomorphic 
$p$-forms on a complex manifold, and $\Omega_{\bar X}^p(\log D)$ will
denote the sheaf of meromorphic $p$-forms with logarithmic poles 
along a divisor.

Given a complex $A^\dt$ with differential $d$, we occasionally denote the 
 cohomology groups by $H^\dt(A^\dt;d)$ when we wish to indicate $d$.

%%%%%%%%%%%%%%%%%%%%%%%%%%%%%%%%%%%%%%%%%%%%%%%%%%%%%%

\head I. Preliminaries on Local Systems.\endhead

\subhead 1. Connections \endsubhead

Let $X$ be a smooth manifold and let ${\cal E}^*_X$ be the sheaf of $\Ci$ $\C$-valued forms
 on it.
More generally, if $V$ is a $\Ci$ complex vector bundle let ${\cal E}^p_X(V)$ be the sheaf of 
$\Ci$ $V$-valued $p$-forms. 
 A connection on $V$ is a $\C$-linear operator
$$\del \colon  {\cal E}^0(V) \to {\cal E}^1(V)$$
satisfying the Leibnitz rule: $\del (fv) = f\del(v)+df \otimes
v$. $\del$ can be extended to a graded 
derivation on ${\cal E}^*(V)$. It can be shown  that $\del^2$ is multiplication by a
fixed $2$-form called the curvature, if it vanishes $\del$ is said to be flat or integrable. For
example, let $V$ be the trivial line bundle, any $1$-form $\theta$ determines a connection
$d-\theta$ on $V$ with curvature $-d\theta$. When
$\del$ is integrable, the differential equation $\del v =0, v(x)=v_0$ has a unique
solution in any contractible neighborhood of $x\in X$. Global
 solutions tend to be multivalued, and the multivaluedness is measured by the monodromy.

Let us spell this out when $V $ is a rank $n$ bundle with a flat 
connection $\del$. 
Choose a Leray covering
$\{U_i\}$ of $X$, that is an open covering such that intersections of a 
finite number of the given sets are either
empty or contractible.  We can find a basis 
 of nonzero sections $\lambda_i^{1},\ldots \lambda_i^{n}$ of
$V|_{U_i}$ such that $\del \lambda_i^{\dt} = 0$. Let $\rho_{ij}$ 
be the change of basis matrix relating $\lambda_{i}|_{U_{ij}}$ and 
$\lambda_{j}|_{U_{ij}}$ on the intersection $U_{ij} = U_i \cap U_j$.
 Then
 $\{\rho_{ij}\}$ is a constant \v Cech $1$-cocycle
with values in $GL_n(\C)$. 
We will denote the corresponding \v Cech cohomology class by 
$$\mu(V,\del) \in H^1(X, Gl_n(\C)) \cong
 Hom(\pi_1(X), Gl_n(\C))/\roman{ conjugacy}.$$

A straight forward argument shows that:

\proclaim {Lemma 1.1} Let ${\cal E}_{cl}^1(X)$ denote the space of closed $1$-forms and let $FC$ denote the set
of flat $\Ci$ connections on the trivial line bundle, then the diagram
$$\matrix
{\cal E}_{cl}^1(X) &{\mathop {\longrightarrow} \limits^{\theta \mapsto d-\theta}}  & FC \\
   [\,\,] \downarrow     &           &  \mu \downarrow\\
 H^1(X,\C)          &{\mathop {\longrightarrow} \limits^{\exp}}  &
H^1(X,\C^*)\\
\endmatrix$$
commutes.
\endproclaim

\proclaim {Variant} The trivial bundle can be replaced by any flat line bundle $(L, \del)$.  
Then the corresponding diagram (with $\del$ replacing $d$, 
the set of flat connections on $L$ replacing $FC$, and 
$\mu(L,\del)\cdot \exp$
replacing $\exp$) commutes.
\endproclaim

We will be primarily interested in the  holomorphic version of all of this. 
So from now on,
 we will assume that $\bar  X$ is a compact complex  manifold,
 $D = \cup \,D_i\subset \bar  X$ a 
divisor with normal crossings and
 $X = \bar  X - D$. Let $j\colon  X \hookrightarrow \bar  X$ denote the
inclusion.  
Let  $V$ be a holomorphic
  vector bundle on $X$. We will usually use the same symbol for the
sheaf of holomorphic sections, 
 although  sometimes we will  write $O_X(V)$ to avoid confusion.
A $\Ci$ connection $\del$ on $V$ is called holomorphic if the derivative,  with respect to $\del$,  
of any local holomorphic section is again holomorphic. This is equivalent to requiring that the
 $(0,1)$
 part of the $\del$ is just $\bar \partial$. 
 In particular there is at most one holomorphic structure on $V$ for
which a given connection is holomorphic.  
Given a pair $(M, \del)$ consisting of a  $\Ci$ vector bundle
and a flat connection, the $(0,1)$ part is automatically integrable 
and hence defines a compatible holomorphic
structure on $M$.

In general, there are an infinite number of holomorphic extensions
 $\bar V$ of $V$ to $\bar X$.
If $V$ carries a flat holomorphic connection $\del$ then there is a unique extension $\bar V$,
 which we
will call the Deligne extension \cite{D1, II 5}, such that $\del$ extends
to a  meromorphic connection
$$\bar \del \colon  \bar V \to \OmD 1 \otimes \bar V$$
with residues having eigenvalues with 
real part in $[0, 1)$. 

Suppose that $(\bar V, \bar \del)$ is a pair as above but without any restrictions on the residues. Then $\bar \del$ 
 extends to graded derivation of $\OmV$ with square $0$. In particular, we obtain a complex
$$0 \to V \ {\mathop {\longrightarrow} \limits^{\bar \del}} \OmD 1 \otimes V 
{\mathop {\longrightarrow} \limits^{\bar \del}} \dots  . $$
Suppose that $\cV = ker  \del (V \to \Om 1 \otimes V)$ is the locally constant sheaf of flat sections on $X$
then  the holomorphic Poincar\'e lemma essentially implies that the restriction of  $\OmV$ to $X$ is quasi-isomorphic
to $\cV$.  Thus we obtain a morphism (in the derived category of bounded 
below complexes):
$$\OmV \to {\R} j_*\cV .$$
\smallskip

\proclaim {Proposition 1.2}  (Deligne \cite{D1, II 6.10}) If none of  the 
eigenvalues of the residues of $\bar \del$ are positive integers then
the above map is an isomorphism, consequently 
$$ H^i(X, \cV) \cong {\HH}^i(\bar X, \OmV) $$
for all i. \endproclaim

When the above conditions are satisfied then as usual we obtain a Hodge to De Rham spectral sequence:
$$ E_1^{pq} = H^q(\bar X, \OmD p \otimes \bar V) \quad \Rightarrow \quad H^{p+q}(X, \cV).$$

A unitary flat vector bundle on $X$ is a pair $(V,\del)$ whose
monodromy representation lies in
$$im[H^1(X,U_n) \to H^1(X, GL_n(\C))].$$

 \proclaim {Theorem 1.3} (Timmerscheidt \cite{T}) If $\bar X$ is compact K\"ahler and 
  $(\bar V, \bar \del )$  the Deligne extension of   unitary flat bundle
$V$ then the above spectral
sequence degenerates at $E_1$.
\endproclaim 

Saito \cite{Sa1,Sa2} has proved much stronger results using variations of
Hodge structures around the same time. Further discussion of all of this will be given in chapter IV.

%************************************************************
\subhead 2. Cohomology Support Loci \endsubhead

In this section, we will study the cohomology support loci which are 
sets of representations of the fundamental group
such that the associated locally constant sheaves
 have nonvanishing cohomological cohomology. We will show that
 these constraints will be 
expressible as algebraic conditions on the coefficients of
representation,  so that
these sets will turn out to be Zariski closed.

First, let us recall the basic facts from invariant theory that 
will be needed.

\proclaim {Theorem 2.1} \cite{MF, pp 27-29}. If $G$ is a reductive group acting
algebraically on a complex affine scheme $X = Spec\,R$ of finite type,
 then the algebra of 
invariants $R^{G}$ is finitely generated.
Its spectrum $X//G = Spec\, R^{G}$ is the quotient of
 $X$ by $G$ in the category of complex affine schemes.
The map is $X \to X//G$ is surjective and the image of a closed $G$
invariant set is again closed.
\endproclaim 

The last part is not explicitly stated in \cite{MF}, so we will
indicate its proof. Given a closed invariant set $W_1\subset X$,
let $p$ be a point of the complement of the  image of $W_1$ in
its closure. Let 
 $W_2$ be the fiber over $p$. Applying \cite{MF, p29 cor.1.2} yields an 
invariant algebraic function $f$ which is $0$ on $W_1$ and $1$ on 
$W_2$, but this contradicts the assumption about $p$. Thus $W_1$ must
have closed image.

\proclaim {Remark} The quotient $X//G$ will not coincide (as a set)
with the set of orbits $X/G$ , unless all the orbits are
closed.
\endproclaim

The above result can be applied to construct the
 so called character variety \cite{LM, S3}.

\proclaim {Theorem 2.2} Let $\Gamma$ be a finitely generated group.
Then the functor $Y \mapsto Hom(\Gamma, Gl_n(H^0(O_Y)))$ is representable
by  an affine scheme $Rep$ of finite type  over $Spec\C$.
The group $Gl_n(\C)$ acts algebraically on $Rep$ by conjugation, and
the complex points of the quotient $H^1(\Gamma, Gl_n(\C))_{ss} =
Rep//Gl_n(\C)$  correspond
to isomorphism classes of semisimple representations of $\Gamma$ into
$Gl_n(\C)$.
\endproclaim

We  will need a description of the character variety from a  \v Cech point
of view. 
Let $X$ be a connected topological space with a finite Leray cover 
${\cal U}=\{U_0, U_1,\ldots U_N\}$.
  The proof of following
result will be given elsewhere.

\proclaim {Theorem 2.3} The set  
 $Z = Z^1({\cal U}, GL_n(\C))$ of constant \v Cech $1$-cocycles has the 
the structure of (the complex points of) an affine scheme in a
 natural fashion. The group
$G = Gl_n(\C)^{N+1}$ acts algebraically on $Z$ such that 
two cocycles are cohomologous if and only if they lie in the same
$G$ orbit. There is an isomorphism
$$H^1(X,GL_n(\C))_{ss} \cong  Z//G.$$
\endproclaim

Given a cocycle
 $ \rho\in Z^1({\cal U},GL_n(\C))$ , define the sheaf
$$(\C^n)_{\rho} = ker[ \prod_c\,j_{c*}\C^n_{U_c} 
{\mathop {\longrightarrow} \limits^{R}} \prod_{a<b}\,j_{ab*}\C^n_{U_{ab}}]$$
where $j_a\colon U_a \hookrightarrow X$ etcetera are the inclusions and
$R((x_c))_{ab} = x_a|_{U_{ab}} - \rho_{ab}x_b|_{U_{ab}}.$
This is a locally constant sheaf of $n$ dimensional vector spaces.
The sheaves $(\C^n)_\rho$ and $(\C^n)_{\rho'}$ are isomorphic if and
only if $\rho$ and $\rho'$ are cohomologous. Furthermore any 
rank $n$ locally constant sheaf of vector spaces is isomorphic to one
of these special sheaves.
Thus $H^1(X,GL_n(\C))$ (respectively $H^1(X, GL_n(\C))_{ss}$) can be viewed
as the parameter space for isomorphism classes of (semisimple)
rank $n$ locally constant sheaves of vector spaces on $X$. 
 
Let $R$ denote the coordinate ring of the affine variety $Z^1({\cal 
U}, GL_n(\C))$. Each cocycle $\rho$ determines a $\C$ algebra 
homomorphism $ev_\rho\colon R \to \C$. 
There is a universal cocycle $\{{\bold u}_{ij}\}\in Z^1({\cal U},GL_n(R))$ 
which satisfies $ev_\rho({\bold u}_{ij}) = \rho_{ij}$.
If $M$ is an $R$-module, then we will
 write $ev_\rho(M)$ for $M\otimes_R \C$  with respect to 
$ev_\rho\colon R \to \C$.

\proclaim {Proposition 2.4} There exists a finite complex of finitely
 generated free $R$-modules
$(V^\bullet, \partial)$ such that
$$ H^i(ev_\rho(V^\bullet), ev_\rho(\partial)) \cong H^i(X, (\C^n)_{\rho}). $$ 
for all $\rho$.
\endproclaim
 
\pf We define
$$ V_{i_0,\dots i_m} =  \cases R^n &\text{if $U_{i_0,\dots i_m} \not= \emptyset$} \cr
                             0 &\text{otherwise.}\cr \endcases$$
Let 
$$V^n =  \bigoplus_{i_0<\dots<i_{m-1}}\, V_{i_0,\dots i_{m-1}}$$
with the differential of $v =(v_{i_0,\dots i_{m-1}}) \in V$ given  by
$$(\partial(v))_{i_0,\dots i_m} =
{\bold u}_{i_0i_1}v_{i_1,\dots i_m} + \sum_{k=1}^m\, (-1)^kv_{i_0,\dots \hat i_k,\dots i_m}.$$
From the definition of $(\C^n)_{\rho}$ it is clear that there are canonical isomorphisms 
$\C^n \cong H^0(U_i,\C_{ \rho}^n)$.
Thus we get a chain of isomorphisms (when $U_{i_0,\dots i_m} \not= \emptyset$),
$$ev_\rho(V_{i_0,\dots i_m}) = \C^n \cong H^0(U_{i_0},\C_{ \rho}^n) 
 \cong H^0(U_{i_0,\dots i_m},\C_{ \rho}^n)$$
where the last map is just restriction. These maps define an isomorphism between 
 $(ev_\rho(V^\bullet),ev_\rho(\partial))$ and the  \v Cech complex 
$C^\bullet({\cal U}, \C_{\rho})$.

\bigskip

An element  $\rho \in H^1(X,GL_n(\C))_{ss}$ determines an isomorphism class 
of locally constant sheaves with semisimple monodromy. We will denote 
a representative of this class by $\C_\rho$. Choose integers $k,m,n$. 
We define the cohomology 
support locus:
$$\Sigma_m^k(X, GL_n(\C)) 
= \{ \rho \in H^1(X, GL_n(\C))_{ss} 
| \, dim \, H^k(X, (\C^n)_{\rho})  \ge m \}. $$
We write $\Sigma_m^k(X)$ for $\Sigma_m^k(X,\C^*)$.

\proclaim {Corollary 2.5} The set 
$\Sigma_m^k(X, GL_n(\C)) $
is Zariski closed in $H^1(X,GL_n(\C))_{ss}$.
\endproclaim 

\pf  As a first step, let us show that the set
$$ \tilde \Sigma =  
\{\tilde \rho \in Z^1({\cal U}, GL_n(\C)) 
| \, dim \, H^k(X, \C_{\tilde \rho}) \ge m \} $$
is Zariski closed.  From the above proposition,  
we conclude that $\tilde \rho$ is in the above set if
 and only if 
$$ rank( ev_{\tilde\rho}(\partial^{i-1})) + rank(ev_{\tilde\rho}(\partial^i )) 
\le rank \, V^i - m.$$
Thus the $\tilde \Sigma$ can be described in terms of the vanishing of 
 certain minors of
 $ev_{\tilde\rho}(\partial^{i-1})$ and $ev_{\tilde\rho}(\partial^i )$,
 and is therefore Zariski closed. This set is also 
 $GL_n(\C)^{\{0\ldots N\}}$ invariant, because $H^k(X,\C_{\tilde 
 \rho})$ depends only on the isomorphism class of $\C_{\tilde\rho}$.
 Therefore the image $\tilde\Sigma$ is Zariski closed.

%*****************************************************************

\head II. Exponentials of Hodge structures.\endhead

Suppose that that $X$ is the complement of divisor in a compact 
K\"ahler manifold. Then $H^1(X,\Z)$ carries a special kind of mixed
Hodge structure that we will refer to as a $1$-Hodge structure.
The principle aim of this chapter is to characterize subsets of
$H^1(X,\C^*) = H^1(X,\Z)\otimes\C^*$
which arise as images of  sub $1$-Hodge structures
 under the exponential map $\exp\colon H^1(X,\C)\to H^1(X,\C^*)$. This 
 result will be 
 applied to cohomology support loci in the last chapter.

\subhead 1. 1-Hodge Structures \endsubhead

A {\it $1$-Hodge structure} is the same thing as a mixed Hodge structure of type 
$\{ (0,1),(1,0),(1,1) \}$ \cite{D}. It consists of:

1) a finitely generated abelian group $\La$.

2) a subspace $W = W_1 \subseteq \La_\Q$.

3) a subspace $F = F^1 \subseteq \La_\C$

\noindent such that 
$W_\C = (W_\C \cap F) \oplus (W_\C \cap \bar F)$
and 
$ \La_\C = W_\C + F.$
\bigskip

The following is elementary.

\proclaim {Lemma 1.1} There are bijections between the following sets of data:\newline
\indent A) $1$-Hodge structures, \newline
\indent B) finitely generated abelian groups $\La$ together with a 
subspace $W \subseteq \La_\Q$ and
a bigrading $\La_\C = H^{01} \oplus H^{10}\oplus H^{11}$ satisfying $\bar H^{01} = H^{10}$,
$\bar H^{11} = H^{11}$ and $W_\C = H^{01} \oplus H^{10}$, \newline
\indent C) finitely generated abelian groups $\La$ together with a subspace $W \subseteq \La_\Q$, a 
complex structure $J$ on $W_\R$ and a complementary subspace $H_\R^{11}$ to $W_\R \subseteq
\La_\R$, \newline
given by:\newline
\indent (A $\to$ B): $H^{10} = F \cap W_\C,\quad H^{01} = \bar F \cap W_\C, \quad
H^{11} = F \cap \bar F$ \newline
\indent (B $\to$ C): $J$ corresponds to multiplication by $i$ under the $\R$-isomorphism
$$ W_\R \to W_\C \to H^{10}$$
and $H_\R^{11} = H^{11} \cap \La_\R$. \newline
\indent (C $\to$ A): $F = (+i\,\, \roman{ eigenspace \,\, of}\, J\otimes \C) + H_\R^{11}\otimes \C.$
 \endproclaim

With the obvious notion of morphism (homomorphisms of abelian groups preserving all the structure)
the above sets become categories and the bijections induce isomorphisms between them. 
We shall be 
primarily  interested in the case where the underlying group $\La$ is torsion free, and this will
assumed from now on (unless otherwise indicated).
\smallskip

\proclaim { Example 1.2} (Hodge-Deligne \cite{D}) Let $\bar X$ be a compact K\"ahler manifold and $D \subset \bar X$
a divisor with normal crossings. Let $X = \bar X - D$ then $\La =
H^1(X, \Z)$ carries a canonical $1$-Hodge
structure.\endproclaim
\smallskip

By abuse of notation, we will say that a subspace $K \subseteq \La_\Q$ is a sub $1$-Hodge structure
of $(L,W,F)$ if $(L \cap K, W \cap K, F \cap K_\C)$ is. By virtue of lemma 1, this is equivalent
to requiring that $W \cap K_\R$ is a complex subspace of $W_\R$ and 
$K_\R = (W_\R \cap K_\R) + (H_\R^{11} \cap K_\R).$
For example $W$ is clearly sub $1$-Hodge structure of $\La$.

We will show how to associate a torus to any $1$-Hodge structure. The construction
 generalizes that  of the Jacobian of a curve and is a special case 
  of  the $1$-motive of a mixed Hodge structure \cite{D}. 

Given a $1$-Hodge structure $\La$, set $\La^* = Hom(\La, \Z)$,
$\La_A^* = Hom(\La,A)$ and $F^* = Hom(F,\C)$. We can identify
$\La^*$ with a subgroup of $F^*$ via the composition $\La^* \to \La_\C^* \to F^*$. The map
$\La_\R^* \to F^*$ is injective because $\La_\C = F + \bar F$, therefore 
$\La^*$ is a discrete subgroup.
Therefore $J(\La) = F^*/\La^*$
is a Lie group, and in fact a complex torus since $\La^*$  contains
a basis of $F^*$. $J$ clearly defines a contravariant functor from the category of $1$-Hodge structures
to the category of complex tori and analytic homomorphisms. The extension of $1$-Hodge structures
$0 \to W \to \La \to \La/W \to 0$
induces an extension of complex tori
$0 \to J(\La/W) \to J(\La) \to J(W) \to 0.$
The group $J(W)$ is compact, and is isomorphic as a real Lie group to $W_R^*/(W \cap \La)^*$.
$J(\La/W)$ is affine.

If $X, \bar X$ and $D$ are as (1.2), then
$$J(H^1(X, \Z)) = {H^0(\bar X, \Omega_{\bar X}^1(log D))^* \over H_1(X,\Z)_{ tors.free}}$$
is just the (generalized) Albanese variety $Alb(X)$. After fixing a base point $x_0 \in X$, we
obtain a holomorphic map, the so called Abel-Jacobi map,
$ \alpha \colon  X \to Alb(X)$
defined by
$$ x \mapsto \int_{x_0}^x \> \roman{ (mod \,\, periods)}.$$

%********************************************************

\subhead 2. The subtorus Theorem\endsubhead

Let $\La$ be a finitely generated torsion free abelian group.  We will denote the homomorphism
$$ \exp \otimes 1\colon  \LC \to \LCs$$
simply by $\exp$. 
 $\LCs$ can be endowed with the structure of an algebraic variety in a very natural way,
namely it can be identified with the set of complex points of $spec \, \C[\La^*]$. Its Zariski tangent 
space at $1$ is isomorphic to $\LC$ via the exponential map $\exp$. Thus the tangent space 
of any subvariety at $1$ can be identified with a subspace  of $\LC$.

When $\La$ is a $1$-Hodge structure, we will say that a subset of
$\LCs$ is of {\it exponential Hodge type} if it is a union of a finite
number of sets of the form $\exp(T_\C + r)$ where $T\subseteq \La$ is
a $1$-Hodge structure and $r \in i\La_\R$. We will give a criterion for
this at the end of this section.

\proclaim {Lemma 2.1} There is a natural bijection between the set of affine subtori of $\LCs$ and the
set  of subspaces of $\La_\Q$. \endproclaim

\pf Given a subspace $V \subseteq \La_\Q$, $spec\, \C[(\La \cap V)^*]$ defines an affine subtorus of
$\LCs$. More geometrically this torus is just $\exp(V_\C) \subseteq \LCs$. Conversely given a 
subtorus $X \subseteq \LCs$, $T_{X,1} \cap \La_\Q$ is a subspace. These are inverse processes.
\bigskip

Let 
$\LC = R_1 \oplus R_2 \oplus \ldots R_n$
be a decomposition into real subspaces. A subset of $ S \subseteq \LC$ will be called
 $R_i$-stable with
respect to the above decomposition provided that if $s=r_1 + \ldots r_n \in S$,
with $r_i \in R_i$, then 
$\epsilon r_i + s\in S $ 
for all sufficiently small $ \epsilon \in \R$. $V \subseteq \LCs$ will be called
$R_i$-stable if $\exp^{-1}(V)$ is.

\proclaim {Lemma 2.2} Given $v \in \LC$ and $\epsilon > 0$, the  Zariski closure
of the image of $(-\epsilon, \epsilon)v$ under $\exp$ is an affine subtorus of $\LCs$.
In particular this set is independent of $\epsilon$.
\endproclaim

\pf Let $V(\epsilon) = \exp((-\epsilon,\epsilon)v).$
Clearly
$V(\epsilon)^{-1} = V(\epsilon)$
and
$V(\epsilon/2).V(\epsilon/2) = V(\epsilon).$
By continuity, the same properties hold for the closures (denoted by 
${\cal C}\ldots$). As the 
Zariski topology is noetherian, 
$${\cal C} V(\epsilon/2^n) = {\cal C} V(\epsilon/2^{n+1})$$
for some $n \ge 0$. Thus ${\cal C} V(\epsilon/2^n)$ is connected subgroup,
or in other words an affine subtorus. In particular, as this set is  closed 
under multiplication
it must contain $$V(\epsilon) = V(\epsilon/2^n)^{2^n}.$$

\proclaim {Remark} By the same reasoning, we  see that the Zariski closure of the
image of an open ball in  a real subspace of $\LC$ under $\exp$ is also an affine subtorus. 
\endproclaim

Let us call a decomposition into $\R$-subspaces 
$ \LC = R_1 \oplus R_2 \oplus R_3$
{\it admissible} provided that 
$ \LC = R_1 \oplus \bar R_1 \oplus R_2 \oplus iR_2$
and 
$ R_3 \cap iR_3 = 0.$
For such a decomposition, we have
$$ dim_\R R_1 + dim_\R R_2 = dim_\R R_3 = dim_\C \LC.$$

\proclaim {Lemma 2.3} Given an admissible decomposition as above, let
$S \subseteq T \subseteq \LC$ be complex subspaces such that 
$T = \sum T\cap R_i$,
$T\cap (R_1 + R_2) \subseteq S$ and $S = \bar S$. Then 
$S = T$.
\endproclaim

\pf As $$ T \supseteq (T\cap R_3) \oplus i(T \cap R_3), $$ we obtain
$ dim_\R T /2 \ge dim_\R T\cap R_3 . $
Therefore, since $T = \sum T\cap R_i$, we obtain
$ dim_\R T /2 \le dim_\R T\cap R_1 + dim_\R T \cap R_2. $
On the other hand the inclusion
$$ (T \cap R_1) \oplus \overline{(T \cap R_1)} \oplus (T\cap R_2) \oplus i(T\cap R_2) 
 \subseteq S $$
implies an inequality
$ dim_\R S /2 \ge dim_\R T\cap R_1 + dim_\R T \cap R_2. $
Combining these shows that $dim_\R S \ge dim_\R T$ which forces equality.

\proclaim {Theorem 2.4} Suppose that $\LC = R_1 \oplus R_2 \oplus R_3$ is
an admissible decomposition. 
 Let $V \subseteq \LCs$ be a Zariski
closed  subset which is $R_1$ and $R_2$ stable. Then any irreducible component of $V$ is of
the form $\exp(T + r)$ where $T$ is a rationally defined linear subspace of $\LC$ such that
$T = \sum T\cap R_i$, and $r\in R_3$. In particular, each component is a translate of 
an affine subtorus. 
\endproclaim

\pf Let $V'$ be an irreducible component of $V$. Choose a smooth point $u \in V'$ 
not lying on any other component, and let $U = u^{-1}V'$. Let $T \subseteq \LC$ be the
tangent space to $U$ at $1$. Given a vector $v\in T$, decompose it as $v_1 + v_2 +v_3$, where
$v_i \in R_i$. Choose a $C^\infty$ curve
$\gamma \colon  (-\epsilon, \epsilon) \to \LC$ such that
$\gamma(0) = 0$, $\gamma'(0)=v$ and the image of $\exp\circ\gamma$ lies in $U$. We can decompose
$\gamma(t) = \gamma_1(t) + \gamma_2(t) + \gamma_3(t)$ with respect $R_1,R_2$ and $R_3$. By assumption 
for any $|t| <
\epsilon$,
 $\exp(s_1\gamma_1(t) + s_2\gamma_2(t) + \gamma_3(t)) \in U$ 
for all $s_i$ close to $1$, and therefore  for any $s_i \in \R$ by lemma 2.2.
This implies that $v_i \in T$, in other words that $T$ is compatible with the decomposition 
into the sum of $R_i$. 
 Furthermore, after taking $s_1 = s_2 = 1/t$ above we obtain
$$\exp(v_1 + v_2) = \lim_{t\to 0} \, \exp(\gamma_1(t)/t + \gamma_2(t)/t + \gamma_3(t)) \in U.$$
That is $\exp(T\cap (R_1 + R_2))$ is contained in $U$, and its Zariski closure $S'$ is an affine 
subtorus (by the remark following 2.2).
 The tangent space $S$ to $S'$ at $1$ is a complex subspace of $T$ satisfying the 
conditions of lemma 2.3. Therefore $S = T$, and consequently $U$ is the  subtorus $\exp(T)$.

If we, write $u = \exp(r_1 + r_2 + r_3)$ with $r_i \in R_i$, then from the preceding 
discussion
$V' = \exp(T + r_3)$.
\bigskip

Suppose that $\La$ is a $1$-Hodge structure then 
$$\LC = H^{10} \oplus H^{11}_\R \oplus i\La_\R$$
is an admissible decomposition.

\proclaim {Corollary 2.5} Let $\La$ be a $1$-Hodge structure and let $V \subseteq \LCs$
be a Zariski closed subset which is $H^{10}$ and $H^{11}_\R$
stable. Then $V$ is of exponential Hodge type.
\endproclaim

\pf By the theorem any component of $V$ is of the form $\exp(T+r)$
with $r$ imaginary and $T$ a rationally defined subspace satisfying
$$T = (H^{10}\cap T) \oplus (H^{11}_\R \cap T) \oplus iT_\R.$$
Therefore the space 
$$ S = (H^{10}\cap T) + (H^{01}\cap T) + (H^{11} \cap T)$$
satisfies the conditions of lemma 2.3, and so $T=S$ which implies  that
$T$ is a sub 1-Hodge structure.

%%%%

\head III. $\bC$-Mixed Hodge Theory. \endhead

In this chapter we introduce the concept of a $\C$-mixed
Hodge structure, which arises from the usual notion by forgetting the
underlying real structure. A detailed development of the basic theory will be 
presented elsewhere. However the essential features of the theory (to 
the extent required here) will be summarized in the first two sections.
The last section contains the key technical results
that will be needed in the next chapter.

\subhead 1. $\bC$-Mixed Hodge structures\endsubhead

 We follow the standard convention that a filtration with superscript is
decreasing, otherwise it is increasing. Any formula involving decreasing filtrations 
can be converted to one involving an increasing filtration by setting $W_n = W^{-n}$.
A $\C${\it-mixed Hodge structure} is a finite dimensional vector space $H$ together
with  three filtrations $(W_\dt, F^\dt,  C^\dt)$ which  are opposed in the
sense of Deligne \cite{D 1.2}, i.e. 
$$Gr_F^p Gr_C^q Gr^W_n H = 0$$
when $p + q \not= n$. The terminology is a little nonstandard but seems 
justified as
any  mixed Hodge structure can be regarded as $\C$-mixed Hodge structure after 
forgetting $\R$-structure  and setting $C =\bar F$. The set of $\C$-mixed Hodge 
structures becomes a category in the obvious way: morphisms are $\C$-linear maps
 preserving all the filtrations. Deligne \cite{D, 1.2.10} shows:

\proclaim {Theorem 1.1} The category of $\C$-mixed Hodge structures is abelian,
 and the filtrations $F$, $C$ and $W$ are strictly preserved by morphisms i.e. the
functors $Gr_F^p,...$ are exact. 
\endproclaim

We say that $H$ is pure of weight $n$
if $Gr^W_nH = H$; this equivalent to statement that  the filtrations $F^\dt$ and 
$C^\dt$ are $n$-opposed. An immediate corollary of the theorem is that a morphism 
between pure objects of different weights must vanish. Finally observe that
if $(F^\dt, C^\dt)$ is $n$-opposed then $(F^{\dt + a}, C^{\dt + b})$ is 
$n-a-b$-opposed

%*****************************************************

\subhead 2. $\bC$-Hodge Complexes \endsubhead

We define a $\C${\it-Hodge complex} to be a quadruple $(A^\dt,F^\dt, C^\dt, W^\dt )$
consisting of   a bounded below complex of
$\C$-vector spaces $A^\dt$ with finite dimensional  cohomology, together
with  three biregular filtrations on it such that:

a) For each $p$, the filtrations on $Gr^p_WA$ induced by $F$ and $C$ are
strictly preserved by the differentials.

b) For each $p,q$, $F$ and $C$ induce $q$-opposed filtrations on
$E_1^{pq}(A, W) = H^{p+q}(Gr_W^pA^\dt)$.

For example the De Rham complex $A^\dt$ of a compact K\"ahler manifold with Hodge
filtration $F$, $C = \bar F$ and $Gr^0_WA =A$ is $\C$-Hodge complex.

The notion of a $\C$-Hodge complex will  be globalized.
 If $S$ is a sheaf on a topological space $X$,
Godement has constructed a canonical flasque resolution ${\cal G}^\dt(S)$ 
\cite{G}. For example
${\cal G}^0(S)$ is the product of stalks of $S$ viewed as sky-scraper sheaves, ${\cal G}^1(S)$
is ${\cal G}^0$ applied to the cokernel of $S \to {\cal G}^0(S)$ and so on. If $S^\dt$ is a 
complex of sheaves  set ${\cal G}^\dt (S^\dt)$ to the associated single 
complex. Then

$\Gamma({\cal G}^\dt (S^\dt))$ is a canonical representative of 
${\R}\Gamma S^\dt$. If 
$(S^\dt, F^\dt)$ is a
 bounded below complex of sheaves with a biregular filtration then 
$({\cal G}^\dt(S^\dt),{\cal G}^\dt(F^\dt))$
is a filtered acyclic resolution (see \cite{D 1.4} for the definition). Thus
$$(\Gamma({\cal G}(S)),\Gamma({\cal G}(F))$$
represents the direct image of $(S,F)$ in the sense of filtered derived categories.
 This can be iterated in the 
presence of several filtrations to a obtain a multifiltered acyclic resolution, and this
leads to a concrete representation for direct images for multifiltered derived categories.

Call a trifiltered complex of sheaves $(A,F,C,W)$ a 
{\it cohomological }$\C${\it-Hodge complex}
provided that $A$ is a bounded below complex of sheaves of $\C$-vector spaces, the
filtrations are biregular and $(\Gamma({\cal G}(A)),\Gamma({\cal G}(F)),\dots)$ is 
$\C$-Hodge complex.

The proofs of the following will be given elsewhere. The first result 
is an analogue of scholium 8.1.9 of \cite{D} and the second is a slight refinement.
					
\proclaim {Proposition 2.1} 
If $(A,F,C,W)$ is a (cohomological) $\C$-Hodge complex then the
spectral sequence
$$E_1^{pq}(A,W) = { H}^{p+q}(Gr_W^pA) \Rightarrow { H}^{p+q}(A)$$
$$({resp.}\> E_1^{pq}(A,W) = {\HH}^{p+q}(Gr_W^pA) \Rightarrow {\HH}^{p+q}(A))$$
degenerates at $E_2$.
\endproclaim

\proclaim {Proposition 2.2} 
Let $(A,F,C,W)$ be a (cohomological) $\C$-Hodge complex 
and $(D,V)$ another
bounded below biregularly filtered complex (of sheaves) of $\C$-vector spaces. Suppose that 
there is a morphism $(F^k,F^k \cap W) \to (D,V)$ such that the sequence
 $$ 0 \to F^{k+1}Gr_W^\dt \to F^kG_W^\dt \to Gr_V^\dt D \to 0$$
is exact. Then $E_2(D,V) = E_\infty(D,V)$ and 
${H}^*(Gr_F^kA) \cong { H}^*(D)$
(resp. ${\HH}^*(Gr_F^kA) \cong {\HH}^*(D)$) noncanonically.
\endproclaim

%**************

\subhead 3. Key Constructions\endsubhead

In this section, we give the basic construction used later on.
Let $R$ be a commutative differential graded algebra over $\C$.
This means that $R$ is a graded vector space with an associative
multiplication which commutes up to sign, and a degree one graded
derivation $d$ satisfying $d^2=0$.
We will assume that $R$ also possesses 3 multiplicative filtrations $F_R,C_R,W_R$.
Suppose that $A$ is a differential graded $R$-module filtered by
a family of three differential graded submodules $F^\dt,C^\dt,W^\dt$ such that
$F_R^pF^i \subseteq F^{p+i}$ etcetera.
Furthermore, we will assume that $(A,F,C,W)$ is a $\C$-Hodge complex.
Our construction depends on a fixed choice of integers $a,b, N$, with N positive,
and two elements 
$\theta \in F_R^a\cap C_R^0\cap W^{-a-b+1}R^1$
and $\psi \in  W_R^{-a-b+1}R^0$ such that $d\theta =0$ and $\theta - d\psi \in C^b$.
Given all this data, we construct  a tri-filtered complex $Cons(A,\theta,\psi)$.

We will first treat the special case where $\psi = 0$, as it is easier and it
motivates the general construction. In this case $\theta \in 
 F^a\cap C^b\cap W^{-a-b+1}R^1$.  As a graded module $Cons(A,\theta,0)$ is
given by
$B^n = (A^n)^{\oplus N}$
  The standard relations
$\theta \wedge \theta = 0$
and
$ d(\theta \wedge a) = -\theta \wedge da$
imply that
$\delta (a_0, a_1, \ldots a_{N-1}) = (da_0, da_1 + \theta \wedge a_0, \ldots
da_{N-1} + \theta \wedge a_{N-2}) $
is a differential on $B$. In other words, 
$(B,\delta)$ is the single complex associated
to the double complex:

$$\matrix   &                             & \ldots& &   &      \cr
          A^0&       \lblarr{\theta\wedge} & A^1& \lblarr{\theta\wedge} &\ldots& A^N \cr
          \downarrow d&                   & \downarrow d& &  & \downarrow d\cr
          A^1&       \lblarr{\theta\wedge} & A^2& \lblarr{\theta\wedge} &\ldots& A^{N+1} \cr
             &                             & \ldots& &   &      \cr 
             \endmatrix $$
or more succinctly
$$B = s(A \lblarr {\theta\wedge} A[1] \lblarr {\theta\wedge} \ldots A[N]). $$
Filter $B$ by
$${\cal F}^p B = s(F^pA \to F^{p+a}A[1] \to \ldots F^{p+Na}A[N]) $$
$${\cal C}^pB = s(C^pA \to C^{p+b}A[1] \to \ldots C^{p+Nb}A[N])$$
and
$${\cal W}^pB = s(W^pA \to W^{p-(a+b)}A[1] \to \ldots W^{p-N(a+b)}A[N]).$$
 In other words, 
$(a_0,\ldots, a_{N-1}) \in \cal F$ if and only if $a_i \in F^{p+ia}$ for $i$,
etcetera. $Cons(A,\theta,0)$ consists of $B$ with these filtrations.

\proclaim {Proposition 3.1}  $Cons(A,\theta,0)$ is a $\C$-Hodge complex.
\endproclaim
		
\noindent {\bf Proof.} By construction
$$Gr_{\cal W}^p B = s(Gr_W^p \lblarr 0 Gr_W^{p-(a+b)} \lblarr 0 \ldots)$$
$$ = \bigoplus_{k=0}^N \, Gr_W^{p-k(a+b)}A $$
$${\cal F}^iGr_{\cal W}^pB = \bigoplus_{k=0}^N \, F^{i+ka}Gr_W^{p-k(a+b)}A $$
$${\cal C}^iGr_{\cal W}^pB = \bigoplus_{k=0}^N \, C^{i+kb}Gr_W^{p-k(a+b)}A $$
Thus $\cal F$ and $\cal C$ are strictly preserved by differentials of $Gr_{\cal W}^p B$.
Furthermore
$$(H^{p+q}(Gr_{\cal W}^pB),{\cal F}^\dt,{\cal C}^\dt)
\cong \bigoplus_{k=0}^N \,(H^{p+q}(Gr_W^{p-k(a+b)}A),F^{\dt +ka},
C^{\dt + kb} )$$
is $q$-opposed.

\bigskip

We now consider the general case where $\psi$ may be nonzero.
We can form a complex $B$ with filtrations ${\cal F}, {\cal W},{\cal C}$
as before, however ${\cal C}^p$ is no longer a subcomplex.
To correct this we define a new filtration $\tilde {\cal C}$
by $(a_0,\ldots, a_{N-1}) \in {\tilde {\cal C}}^p$ if and only if
for each $i$,
$$ a_i + \psi a_{i-1} + \ldots {\psi^i \over i!}a_0 \in C^{p+ib}.$$
Set 
$Cons(A,\theta, \psi) = (B, {\cal F}, {\tilde {\cal C}}, {\cal W}).$

\proclaim {Lemma 3.2} ${\tilde {\cal C}}^p$ is a subcomplex.
\endproclaim

\pf Suppose $(a_0, a_1,\ldots ) \in {\tilde {\cal C}}^pB^n$ then
$$(da_i + \theta \wedge a_{i-1}) + \psi(da_{i-1} + \theta \wedge a_{i-2})
+\ldots {\psi^i \over i!} da_0 $$
$$ = d(a_i + \psi a_{i-1} + \ldots {\psi^i \over i!} a_0)
+ (\theta  - d\psi)(a_{i-1} + \ldots {\psi^{i-1} \over {i-1}!}a_0) \in
C^{p+ib} $$
because  $\theta - d\psi \in C^b$.

\proclaim {Lemma 3.3}  $ {\tilde {\cal C}}^p Gr_{\cal W}^k B = {\cal
C}^p Gr_{\cal W}^k B$.
\endproclaim

\pf Suppose $(a_0, \ldots a_{N-1}) \in {\cal W}^{k}$ then
$a_i \in W^{k-i(a+b)}$. As $\psi \in W^{-a-b+1}$, we have
$$\psi a_{i-1} + \ldots {\psi^{i-1} \over {i-1}!}a_0 \in W^{k-i(a+b)+1}.$$
Therefore
$$a_i \in C^{p+ib} + W^{k+1-i(a+b)}$$
if and only if 
$$ a_i + \psi a_{i-1} + \ldots {\psi^{i-1} \over {i-1}!}a_0 
\in C^{p+ib} + W^{k+1-i(a+b)}. $$
Thus
$$ {\tilde {\cal C}}^p \cap {\cal W}^k + {\cal W}^{k+1}
= {\cal C}^p \cap {\cal W}^k + {\cal W}^{k+1}, $$
and the lemma follows by dividing both sides by ${\cal W}^{k+1}$.

\proclaim {Proposition 3.4} $(B, {\cal F}, {\tilde {\cal C}}, {\cal W})$ is a 
$\C$-Hodge complex.
\endproclaim

\pf The proof is similar to that of proposition 3.1. The filtrations  $\cal W$ and $\cal F$
satisfy
$$Gr_{\cal W}^p B = \bigoplus_{k=0}^N \, Gr_W^{p-k(a+b)}A $$
$${\cal F}^iGr_{\cal W}^pB = \bigoplus_{k=0}^N \, F^{i+ka}Gr_W^{p-k(a+b)}A $$
By lemma 3.3, $\tilde {\cal C}$ and $\cal C$ induce the same filtration on
$Gr_{\cal W}^\dt B$, therefore
$${\tilde {\cal C}}^iGr_{\cal W}^pB = \bigoplus_{k=0}^N \, C^{i+kb}Gr_W^{p-k(a+b)}A. $$
The remainder of the proof is identical to that of proposition 3.1.
\bigskip

We will now consider the case where
 $R$ is a commutative differential bigraded algebra and $A$ is a differential bigraded $R$-module.
 This means
that $R$ and $A$ are bigraded  with 2 anticommuting
 differentials $d'$ and $d''$, respectively of bidegree $(1,0)$
and $(0,1)$, which satisfy Leibnitz rules.  $R$ becomes
a differential graded algebra and $A$ a differential graded module with total grading
$$R^n = \bigoplus_{p+q=n}\, R^{pq}$$
$$A^n = \bigoplus_{p+q=n}\, A^{pq}$$
and  differential $d = d'+d''$.  
For simplicity
we will assume that $R^{pq}$ and $A^{pq}$ vanishes for all but
finitely many $(p,q)$ and that
$$F_R^k = \bigoplus_{p \ge k} \, R^{pq}$$
and
$$F^k = \bigoplus_{p \ge k} \, A^{pq}$$
Let $\theta \in C^0\cap W^0\cap R^{10}$,
 $\phi \in C^0\cap W^{-1}\cap R^{10}$, and $\psi \in W^{-1}R^{00}$
be elements such that $d\theta=d\phi=0$ and $\phi - d\psi \in C^1$.
Choose $N >>0$.
Then  we can construct $\C$-Hodge complexes
$\calA = Cons(A,\theta, 0)$ with $a=1$ and $b=0$, and
$\calB = Cons(A,\phi, \psi)$ with
$a=b=1$ and the same $N$. 
To make this explicit set 
$$\calA^{pqs} = \calB^{pqs} =A^{p+q,s}$$
with $q = 1,2,\ldots N$.
Then
$$\calA^n= \calB^n= \bigoplus_{p+q+s=n}\, \calA^{pqs}.$$

 The filtrations $\cal F$ and $\cal W$
are given by 
$${\cal F}^k\calA ={\cal F}^kB = \bigoplus_{p \ge k}\, \calA^{pqs}$$
$${\cal W}^k\calA  = \bigoplus W^{k-q} \calA^{pqs}$$
$${\cal W}^kB  = \bigoplus W^{k-2q} \calB^{pqs}$$

Define a map 
$$\sigma \in Hom (\calB^\dt,A^\dt)=Hom( \calA^\dt, A^\dt) = \bigoplus\, Hom(\calA^{pqs}, A^{ij}) $$
as the sum of identity maps  $\calA^{pqs} \to A^{p+q,s}=\calA^{pqs}$. The map 
$\sigma\colon \calA^\dt\to A^\dt$
can be regarded as a map of complexes where $A$ is equipped with the differential $d+\theta$.
Similarly $\sigma\colon  \calB^\dt \to A^\dt$ is a map of complexes when $A$ is given $d+\phi$ as its
differential.
Define  new filtrations
$${\cal V}(b)^k A^{ij} = W^{k-(b+1)i}A^{ij}$$
with $b=0,1$.
Then
$$ \sigma({\cal W}^k\calA^{0qs}) = {\cal V}^k(0)A^{qs}$$
$$ \sigma({\cal W}^k\calA^{pqs}) \subseteq {\cal V}^{k+1}(0)A^{p+q,s}$$
when $p$ is positive.
Therefore the sequence
$$0 \to {\cal F}^1 Gr_{\cal W}^\dt \calA^\dt \to
 {\cal F}^0 Gr_{\cal W}^\dt \calA^\dt \to Gr_{ {\cal  V}(0)} A^\dt \to 0$$
is exact. By a similar argument
$$0 \to {\cal F}^1 Gr_{\cal W}^\dt \calB^\dt \to
 {\cal F}^0 Gr_{\cal W}^\dt \calB^\dt \to Gr_{ {\cal  V}(1)} \calB^\dt \to 0$$
is also exact.
Observe that $Gr_{\cal F}^0\calA^\dt$ is just $A$ with $d''+\theta$ as differential.
Proposition 2.2   implies that 
$$H^\dt (Gr_{\cal F}^0\calA^\dt) =  H^\dt(A^\dt; d''+ \theta)$$
and 
$$H^\dt(A^\dt; d+\theta)$$
are isomorphic. Furthermore the spectral sequences abutting to either group with
$$E_1^{pq} = H^{p+q}(Gr_{{\cal V}(0)}^p \calA)$$
degenerates at $E_2$. Similarly
$$ H^\dt(A^\dt; d'' +\phi) \cong H^\dt(A^\dt; d+\phi),$$
and that the spectral sequence associated to ${\cal V}(1)$ degenerates at $E_2$.

\proclaim {Proposition 3.5} If $H^i(A^\dt; d+\theta +\phi) \not= 0$
then $H^i(A^\dt; d+\phi) \not=0$.
\endproclaim

\pf  Let $\delta = d+ \theta + \phi$ or $\delta =  \phi$ then there
are spectral sequences
$$E_1^{pq}(\delta) = H^{p+q}(Gr_{{\cal V}(1)}^pA^\dt) \Rightarrow H^{p+q}(A^\dt,\delta).$$
The $E_1$ differential is just the connecting map associated to
$$ 0 \to Gr_{{\cal V}(1)}^{p+1} (A^\dt,\delta) \to 
{\cal V}(1)^p(A^\dt,\delta)/{\cal V}(1)^{p+2}(A^\dt,\delta)
\to Gr_{{\cal V}(1)}^{p+1}(A^\dt,\delta) \to 0 . $$
As $\theta{\cal V}(1)^p \subseteq  {\cal V}(1)^{p+2}$, the exact sequences associated to
$\delta = d + \theta + \phi$ or $\delta = \phi$  coincide. Thus the $E_2$ terms of 
the two spectral 
sequences agree. The hypothesis implies that 
$E_2 ^{p,i-p}( \phi ) = E_2^{p,i-p}(\del +\theta + \phi) \not= 0$
for some $p$. The proposition follows from the equality $E_2(\phi) = E_\infty(\phi)$.
\bigskip

Finally, we will globalize these results. The notation will be almost the same 
as before, but we will repeat it for clarity.
 Let $R^{\dt\dt}$ be a sheaf of commutative differential bigraded
algebras and $A$ a sheaf of differential bigraded $R$-modules with differentials
$d'$ and $d''$ of types $(1,0)$ and $(0,1)$ respectively, and total differential $d=d'+d''$.
 We will assume 
that $R^{pq}$ and $A^{pq}$ vanish for all but finitely many $(p,q)$. Define
$$F_R^k = \bigoplus_{p \ge k} \, R^{pq}$$
and
$$F^k = \bigoplus_{p \ge k} \, A^{pq}.$$
Let $C_R,W_R$  be additional multiplicative filtrations on
$R$. Let $C,W$ be compatible filtrations on $A$ (i.e. $C_R^iC^j\subseteq C^{i+j}$ etc.)
such that $(A,F,C,W)$ is a cohomological $\C$-Hodge complex.

\proclaim {Theorem 3.6}
Let $\theta \in H^0(C^0\cap W^0\cap R^{10})$, $\phi \in H^0(C^0\cap W^{-1}\cap R^{10})$
and $\psi \in H^0(W^{-1}R^{00})$ be elements such that $d\theta =d\phi=0$ and
$\phi-d\psi=0$. Then \hfil \break
\indent 1)
$${\HH}^\dt (A^\dt;d  + \theta) \cong {\HH}^\dt (A^\dt; d''+ \theta), $$
and furthermore the spectral sequences converging to  either hypercohomology group associated
to ${\cal V}(0)$ degenerates at $E_2$. \hfil \break
\indent 2)
$${\HH}^\dt (A^\dt;d  + \phi) \cong {\HH}^\dt (A^\dt d'' +\phi), $$
and furthermore the spectral sequences converging to  either hypercohomology group associated
to ${\cal V}(1)$ degenerates at $E_2$. \hfil \break
\indent 3)
If ${\HH}^i (A^\dt; d  +\theta +\phi) \not= 0 $ then ${\HH}^i (A^\dt; d+\phi) \not= 0$.
\endproclaim

%%%%%%%%%%%%%%%%%%%%%

\head IV. Hodge Decomposition for Higgs Bundles.\endhead

 In this chapter, we will establish the quasi-isomorphisms
 $${\HH}^\dt(\OmV ; \del+\theta) \cong 
{\HH}^\dt(\OmV ; \theta)$$
stated in the introduction. This should be viewed as a strong 
generalization of the classical Hodge decomposition
which corresponds to the case where 
$(V,\del,\theta,D)=(O_X,d,0,\emptyset)$.

\subhead 1. Real Analytic Log Complex\endsubhead

Throughout this chapter, let $\bar X$ denote a compact K\"ahler
manifold, $D \subset \bar X$ a  reduced divisor
with normal crossings, and $X= \bar X - D$.
Let $j \colon  X = \bar X-D \hookrightarrow \bar X$ be the 
inclusion.
The complex $\Omdt$ is filtered by the stupid filtration
$F$ and the weight filtration 
$$W_k\Omega_{\bar X}^i(log D)= \Omega_{\bar X}^{i-k} \wedge \Omega_{\bar X}^k(logD).$$
Navarro Aznar \cite{N, sect. 8} has defined a related complex 
 $A_{\bar X}^\dt(log D)$ consisting 
of  differentials which are real analytic on $X-D$ with certain
allowable  singularities along $D$. 
$A_{\bar X}^\dt(logD)$ is the algebra over the ring of real analytic  functions generated by
 the real analytic differentials and $\log|f|, Re\,df/f, Im\, df/f$ for
any (local) holomorphic function $f$ vanishing along a component of $D$.
The useful feature of this complex is that it stable under complex
conjugation.  It has a bigrading 
$$A_{\bar X}^n(logD) = \bigoplus_{p+q=n}\, A_{\bar X}^{p,q}(logD)$$ 
by $(p,q)$ type, and as usual the  exterior derivative $d$ breaks up into a
sum of operators $\partial + \bar \partial$ of degrees $(1,0)$ and
 $(0,1)$ respectively.
It has a multiplicative weight filtration defined by
$$ W_k A_{\bar X}^{m}(\log D) = im(A_{\bar X}^{k}(\log D) \wedge A_{\bar X}^{m-k} \to 
A_{\bar X}^{m}(\log D) ) $$ 
where $A_{\bar X}^\dt$ is the complex of real analytic forms.
There is a residue map:
$$Res_k\colon  W_k A_{\bar X}^m(log D) \to A_{D^k}^{m-k}$$
where $D^k$ is the disjoint union of $k$-fold intersections of components of
$D$. Suppose $z_i$ is a system of local coordinates such that $D$ is defined by
the vanishing of the product of the first $d$ coordinates. If $I$ is the finite set of integers 
consisting of $i_1 < i_2 < \ldots i_n$, let us write
$$ dz_I/z_I = dz_{i_1}/z_{i_1} \wedge \dots dz_{i_n}/z_{i_n}.$$
 The residue map is determined by the rules: $Res_k$ is additive, if $\omega$ 
is a product of an analytic differential $\alpha$ with $dz_I/z_I$ and
$d\bar z_J/\bar z_J$ with  $I \cap J = \emptyset$ and $card(I \cup J) = k$
then $Res_k(\omega)$ is $\pm \alpha|_{\{\prod z_i \prod z_j = 0\}}$, otherwise
if $\omega$ is sum of terms which not of the previous type then $Res_k(\omega) = 0$.
See \cite{N p45} for the precise rule concerning signs which is given as an explicit
function $\varepsilon(I,J)$.

 Suppose that $D = T + N$, then we can define a partial residue:
$$Res'_k\colon W_k  A_{\bar X}^m(log D) \to A_{T^k}^{m-k}(log N \cap T^k),$$
by modifying the above rule. Suppose $T$ is locally defined by $z_1\dots z_t = 0$.
 If $\omega$ is a product of  differential
forms $\alpha$, $dz_I/z_I$ and $d\bar z_J/\bar z_J$, where $\alpha$ has no
singularities along $T$ and $I,J$ are disjoint subsets of $\{1,\dots t\}$
such that $card(I \cup J) = k$,
then 
$$Res'_k(\omega) = \varepsilon(I,J) \alpha|_{\{\prod z_i \prod z_j = 0\}}.$$

Let $(\bar L, \bar \del)$ be the Deligne extension of a line bundle
with a unitary flat
connection as in section 1 of chapter I . Let $T$ be the union of components $D$ about
 which $\del$ has
trivial monodromy, and let $N = D -T$. The connection $\del$ extends to graded derivation of
$\A {\dt} $ such that $\del^2 = 0$. The $(0,1)$ part of $\del$ is just
the Cauchy-Riemann operator $\bar \partial$, we define $\partial = \del - \bar \partial$.
Thus we obtain a double complex $\A {\dt\dt} $. We define a new weight filtration by
$$ W_k(\A {m} ) = image(A_{\bar X}^{m-k}(\log N) \wedge \A {k} \to \A {m} ) $$
We define a residue map
$$ Res_k\colon  W_k (\A {m} ) \to A_{T^k}^{m-k}(\log N \cap T^k) \otimes {\bar L}|_{T^k}$$
by 
$$ Res_k(\alpha \otimes \lambda) = Res'(\alpha) \otimes \lambda|_{T^k} . $$
We define additional filtrations 
$$F^p \A {\dt} = \bigoplus_{p' \ge p} \, \A {p,q} $$ 

$$ \bar F^q \A {\dt} = \bigoplus_{q' \ge q} \, \A {p,q} $$
The following properties are immediate from the definition:

a) $Res$ commutes with $d$ and with $\omega \mapsto \theta \wedge \omega $ for
any analytic form $\theta$.

b) If we write $F^pW_k  = F^p \cap W_k$ then
$$ Res_k(F^p W_k (\A {\dt} ) )  \subseteq F^{p-k}A_{T^k}^\dt (\log N \cap T^k)
 \otimes {\bar L}[-k] ,$$
and a similar formula holds for $\bar F^p$.

All of this can be carried out for higher rank bundles at the expense
of more notation. Let $(\bar V, \bar \del)$ be 
the Deligne extension of a holomorphic vector bundle with unitary flat 
connection. Let $\cV$ be the corresponding unitary local system on
 $X$. On $D^k- D^{k+1}$ there are unitary flat bundles 
 $(V_k,\del_k)$ corresponding to the  local systems 
 $$  \cV_k = j_{k*}\cV |_{D^k - D^{k+1}}, $$
where $j_k\colon D^k - D^{k+1}\hookrightarrow D^k$ is the inclusion.
Furthermore $V_k$ is a subbundle of the restriction of $\bar V$ to
 $D^k-D^{k+1}$ \cite{T}. 
The double complex $\AV {\dt  \dt} $ can be constructed  with
filtrations $F$ and $\bar F$ as above.
 The definition of the weight filtration is more involved.
The basic idea is that since $\cV$ is unitary, the local monodromy
around $D$ can be diagonalized. Therefore $(\bar V,\bar\del)$ is locally
isomorphic to a direct sum of unitary flat line bundles. There is one and
only one way to choose a filtration
 $$W_k(\AV {\dt} )$$
so as to coincide locally with the direct sum of the $W_k$'s for each line
bundle.  In a similar fashion a residue map
$$ Res_k \colon  W_k(\AV {m} ) \to W_0(A_{D^k}^{m-k}(log D^{k+1})\otimes \bar V_k) $$
can be defined so that the above properties hold.

The complex of holomorphic differentials is a subcomplex of $\A {\dt} $ via
$$ \OmV = ker(\bar \partial \colon  \AV {\dt, 0} \to \AV {\dt, 1} ) .$$
The  filtrations $W$ andf $F$ restrict to this subcomplex. Note that
$F$ restricts to the stupid filtration, and 
$$W_k \OmV = W_k \cap \OmV$$ 
is the weight filtration defined in \cite{T}.

\proclaim {Theorem 1.1} The bifiltered complexes $(\AV {\dt} , W, F)$ and
$(\OmV , W, F)$ are bifiltered quasiisomorphic.
\endproclaim

This was proved by Navarro Aznar \cite{N, 8.8} when $V$ is a trivial line
bundle, and the
general case goes through with only a few modifications which we will
briefly indicate.  The main point is to prove the exactness of
$$ 0 \to W_k(\OmDV p ) \to W_k( \AV {p, 0} ) \lblarr {{\bar \partial }}
 W_k( \AV {p, 1} ) \ldots .$$

First we have to set up the notation.
The problem is local, so we can replace $\bar X$ by small coordinate
neighbourhood of a point $0$.
 We can also
assume that $(\bar V,\bar \del)$ is a sum of  flat line bundles. We may further
reduce to the situation where $\bar V$ is itself a line bundle, which
can be assumed to be trivial (although $\del$ need not be).
We decompose
subdivisors $D = N+T$ as above.  
 Let $z_i$ be a system of local coordinates
such that $N = \{z_1 \ldots z_n = 0 \}$ and $T = \{ z_{n+1} \ldots z_t
= 0 \}$. To avoid a conflict of notation, we will let
$$W'_k A_{\bar X}^{m}(\log D)  =  
im(A_{\bar X}^{k}(\log D) \wedge A_{\bar X}^{m-k} \to 
A_{\bar X}^{m}(\log D) ) $$
$$W_k  A_{\bar X}^{m}(\log D)
 = im(  A_{\bar X}^{k}(\log D) \wedge A_{\bar X}^{m-k}(\log N) \to
 A_{\bar X}^{m}(\log D) ) .$$
Note that $W'$ is denoted by $W$ in \cite{N}. Let us call a differential
form basic,  if it can expressed as
 $$\alpha \wedge \prod_{l} (\log |z_l|)^{k_l}
 {dz_I \over z_I} \wedge {d\bar z_J \over \bar z_J}$$
with $\alpha$  analytic and not divisible by $z_i, \> i\in I$ or 
 $\bar z_j,\> j \in J$.
 Define its $l$th weight of the above form as
$$\cases k_l &\text{ if $l\notin I\cup J$}\cr
         k_l+1 &\text{ if $l \in I\cup J - I\cap J$}\cr
         k_l+2 &\text{ otherwise} \endcases $$
 Define the $T$-weight of a basic form to be the sum of
 $l$th weights with $l \in \{n+1, \ldots t \}$
 By assumption the elements in $A_{\bar X}^{k}(\log D)$ are sums of
finitely many basic forms. Furthermore a form lies in 
$W_k$ if and only if it is a sum of basic forms with  $T$-weight at
most $k$.  A simple calculation shows that if
$\beta$ is basic then $\bar \partial \beta$ is either zero or a sum of
basic forms with the same $i$th weight as
$\beta$.  

Suppose that $ q > 0$ and that  
$$\alpha \in W_k A_{\bar X}^{pq}(\log D ) \subseteq 
 W'_{k+n} A_{\bar X}^{pq}(\log D ) $$
satisfies $\bar \partial \alpha = 0$. Then Navarro Aznar \cite{N} shows
that 
$\alpha = \bar \partial \beta$ with
$\beta \in W'_{k+n} A_{\bar X}^{pq}(\log D ) $. Let $\beta' $ be the
sum of basic components of $\beta$ which lie in $W_k$. Then $\beta''
= \beta - \beta'$ is a sum of basic forms with $T$-weight greater than
$k$.  Therefore  $\bar \partial \beta'' = 0$. Thus 
$  \alpha = \bar \partial \beta'$.

\proclaim {Proposition 1.2} There are quasiisomorphisms induced by the
residue map 
$$ Res \colon    Gr_k^W( \AV {\dt} ) \cong W_0(A_{D^k}^\dt (\log D^{k+1}) 
\otimes \bar V_k) [-k],$$
$$ Res(F^p Gr_k^W( \AV {\dt} ) \cong F^{p-k} W_0 (A_{D^k}^\dt (\log D^{k+1}) 
\otimes \bar V_k) [-k],$$
$$ Res(\bar F^p Gr_k^W( \AV {\dt} ) \cong
 \bar F^{p-k}  W_0(A_{D^k}^\dt (\log D^{k+1}) \otimes \bar V_k) [-k].$$
 \endproclaim 

\pf It suffices to prove the second isomophism, since the first is a special case and
the third follows by conjugation. We will prove the equivalent 
statement
$$Res \colon  Gr_F^p Gr_k^W(\AV ^{\dt} ) \cong Gr_F^{p-k} W_0 A_{D^k}^\dt (\log D^{k+1}) 
\otimes \bar V|_k [-k].$$
By theorem 1.1 it suffices to establish that $Res$ induces an isomorphism
$$ Gr_W^k( \OmDV {p}  )\cong W_0(\Omega_{D^k}^{p-k}(\log D^{k+1})\otimes\bar V_k) ,$$
but this is proved by Timmerscheidt \cite{T, 1.7}.

 Finally, a straight forward modification of the arguments of 
  Timmerscheidt \cite{T, T1} shows that:

\proclaim {Theorem 1.3}  $(\AV {\dt} ,F,\bar F, W)$ is a cohomological
$\C$-Hodge complex
\endproclaim

The idea is as follows. 
By theorem 1.1 and work of Timmerscheidt, the $i$th cohomology
of $W_0(A_{D^k}^\dt (\log D^{k+1}) \otimes \bar V_k)$ is isomorphic to
 the space $H^i_{(2)}(D^k,V_k)$ of harmonic $V$-valued $i$-forms 
satisfying $L^2$ growth conditions
with respect to a suitable K\"ahler metric on $D^k-D^{k+1}$ \cite{Z, sect 3}.
 The  K\"ahler identities imply a decomposition into $(p,q)$ type
$$ H^i_{(2)}(D^k,V_k) = \oplus_{p+q=i} H^{pq}_{(2)}(D^k,V_k)$$
and
$$ H^{pq}_{(2)}(D^k,V_k^*) = \bar H^{qp}_{(2)}(D^k,V_k) .$$
Furthermore, under the above isomorphism
$$H^i(F^mW_0(A_{D^k}^\dt (\log D^{k+1}) \otimes \bar V_k)) \cong
 \oplus_{p \ge m}
H^{p,i-p}_{(2)}(D^k,V_k). $$
As $V$ carries a unitary metric, we can identify $V^*$ with the
complex conjugate of $V$, therefore
$$H^i(\bar F^mW_0(A_{D^k}^\dt (\log D^{k+1}) \otimes \bar V_k)) \cong
 \oplus_{q \ge m}
H^{i-q,q}_{(2)}(D^k,V_k). $$
Thus $F$ and $\bar F$ are strict and induce $i$-opposed
 filtrations. This together with  proposition 1.2 concludes the argument.

%*******************************************

\subhead 2. Main theorem\endsubhead

Let $(\bar V,\bar\del )$ be the Deligne extension of a unitary flat vector 
bundle on $X$. A section 
$$\theta \in H^0(\Omega_{\bar X}^1(log D)\otimes End(\bar V))$$
is called a Higgs field if $\theta\wedge \theta$ viewed as a
section of $\Omega_{\bar X}^2(log D)\otimes End(\bar V)$ vanishes.
Given a Higgs field, the graded sheaf $\OmV$ becomes a complex with  
differential $\xi \mapsto \theta\wedge\xi.$
$\theta$ is closed with respect to the induced unitary connection
on $End(\bar V)$ by theorem 1.3 of chapter I. Therefore 
$\xi \mapsto (\bar \del + \theta\wedge)\xi$
gives a second differential on $\OmV$

\proclaim {Theorem 2.1} If 
$$\theta \in H^0(\Omega_{\bar X}^1\otimes End(\bar V))$$
is a holomorphic Higgs field then 
there are  noncanonical isomorphisms
$${\HH}^*(\OmV ; \theta) \cong {\HH}^*(\OmV ; \bar\del+\theta)$$
\endproclaim 

\pf  This will follow from the results of section 3 of chapter III.
Let $R^{\dt\dt}$ be the subalgebra of 
$A_{\bar X}^{\dt\dt}(log D)\otimes End(\bar V)$
generated by $A_{\bar X}^{\dt\dt}(log D)\otimes Id$ and $1\otimes 
\theta$. This is  a commutative differential bigraded algebra. 
Set 
$A^{\dt\dt} = \AV {\dt\dt} $, $d'= \partial $ and $d''= \bar\partial$.
The first part of theorem of 3.6 of the indicated section implies that
$$ {\HH}^\dt(\AV {\dt} ; \bar\del+\theta) \cong 
{\HH}^{\dt}(\AV {\dt} ; \bar \partial + \theta).$$
Theorem 1.1 of this chapter implies that
$$(\OmV ; \bar\del + \theta) \hookrightarrow (\AV {\dt} ; \bar\del + \theta)$$
and
$$(\OmV ; \theta) \hookrightarrow (\AV {\dt} ; \bar \partial + \theta)$$
are quasi-isomorphisms. So the theorem is proved.

\proclaim {Corollary 2.2} If $t \in \C^*$ then 
$${\HH}^*(\OmV ; \bar\del+\theta) \cong {\HH}^*(\OmV ; \bar\del+ t\theta)$$
\endproclaim

It would be very interesting to allow the Higgs field to have poles.
We will take a small step in this direction which will be needed in
the next section. Let
$$\Phi \in H^{11}(X, \C) = F^1\cap {\bar F}^1H^1(X,\C).$$
Then there exist differential forms $\phi, \phi'\in  F^{1}H^0(\OmD {1} )$ such
$\phi$ and $\bar \phi'$ both represent $\Phi$.

\proclaim  {Lemma 2.3} There exists $\psi \in H^0(W_1A_{\bar X}^{00}(log D))$ such that
$d\psi = \phi - \bar \phi'$.
\endproclaim 

\pf As $H^*(X) = H^*(A_{\bar X}^{\dt}(log D) )$, $\phi - \bar \phi' = d\psi$ 
for some $\psi \in H^0(A_{\bar X}^{00}(log D))$.
Let $z_i$ be local coordinates such that $D$ is defined by $z_1...z_d=0$.
 Then
$\psi$ is a sum of nonzero terms of the form 
$$f_{k_1\ldots k_d}\prod(\log|z_i|)^{k_i}$$
 with $f_{k_1\ldots k_d}$ analytic. The sum of derivatives of these terms lies in 
$W_1$ if and only if
$\sum k_i \le 1$, which implies that $\psi \in W_1$.

\proclaim {Theorem 2.4} If $(\bar L, \bar \del)$ is the Deligne extension of
a unitary flat line bundle on $X$, there is a  noncanonical isomorphism
$${\HH}^*(\OmL ;\phi) \cong  {\HH}^*(\OmL ; \bar\del+ \phi).$$
If $\theta$ is a holomorphic 1-form on $\bar X$ such that
 ${\HH}^i(\OmL, \bar\del + \theta + \phi) \not= 0$ then
${\HH}^i(\OmL, \phi) \not= 0$.
\endproclaim 
 
\pf 
Set $R^{\dt\dt} = A_{\bar X}^{\dt\dt}(log D)$
$A^{\dt\dt} = \A {\dt\dt} $, $d'= \partial $ and $d''= \bar\partial$.
Theorem of 3.6 of chapter III implies that
$$ {\HH}^\dt(\A {\dt} ; \bar\del +\phi) \cong {\HH}^{\dt}(\A {\dt} ; \bar 
\partial + \phi).$$
It follows from theorem 1.1  that
$$(\OmL ; \bar\del + \phi) \hookrightarrow (\A {\dt} ; \bar\del + \phi)$$
and
$$(\OmL ; \phi) \hookrightarrow (\A {\dt} ; \bar \partial + \phi)$$
are quasi-isomorphisms.
This implies the first part of the theorem, and the second
part is similar.

\proclaim {Corollary 2.5} If $t \in \C^*$ then 
$${\HH}^*(\OmL ;\bar\del+\phi) \cong {\HH}^*(\OmL ; \bar\del+ t\phi).$$
\endproclaim 

%%%%%%%%%%%%%%%%%%%%%%%%%%%%%

\head V. Structure of Rank One Loci.\endhead

Let $\bar X$ be a compact K\"ahler manifold, $D \subset \bar X$ a divisor
with normal crossings, and $X= \bar X - D$. A second K\"ahler compactification
of $X$ which dominates $\bar X$ will be called an {\it enlargement} of $\bar X$.

\proclaim {Theorem 1.1} If the mixed Hodge structure on $H^1(X,\Z)$
is pure
then for each pair of integers $k,m \ge0$ $\Sigma_m^k(X)$ is of exponential
Hodge type.\endproclaim 

\pf We use the criterion of chapter II corollary 2.5. Corollary
2.5 of chapter I implies that $\Sigma^k(X)$ is Zariski closed.
If $H^1(X,\Z)$ is pure then either it has weight one (so that
$H^{11}=0$) or  weight two (so that  $H^{10}=0$).
The stability follows from corollaries 2.2 and 2.5 of chapter IV.

\proclaim {Remark 1} A sufficient condition for $H^1(X,\Z)$ to be pure
of weight one is that $X$ has a possibly singular compactification
whose boundary has codimension 2 or more. $H^1(X,\Z)$ is pure
of weight two if $H^1(\bar X,\C) = 0$.
\endproclaim 

\proclaim {Remark 2} By the same sort of argument as above,
 $\Sigma_m^k(X) \cap H^1(\bar X, \C^*)$
is of exponential Hodge type.
\endproclaim

\proclaim {Corollary 1.2} If either $X = \bar X$ or $H^1(\bar X, \C) = 0$ then
for each integer $k \ge 0$, there exists a finite number of unitary characters
$\rho_i \in H^1(X, \C^*)$, and holomorphic maps into complex  tori
$f_i\colon X \to T_i$ such that 
$$\Sigma_m^k(X) = \bigcup_i \rho_if_i^*H^1(T_i, \C^*).$$
\endproclaim 

\pf $\Sigma_m^k(X)$ is expressible as a union of finitely many affine tori
of the form $\exp( \La_i\otimes \C +r_i)$ where $\La_i\subseteq H^1(X,\Z)$ is
a sub 1-Hodge structure, and $r_i\in \La_i\otimes \R$.
 Let $\rho_i = \exp(r_i)$, $T_i = J(\La_i)$ and
$f_i$ be the composition of the  Abel-Jacobi map with the projection
$Alb(X) \to T_i$.
\bigskip

By a curve, we will mean a compact complex curve with finitely many
points removed.
Let us call a holomorphic map $f\colon X \to C$ to a smooth curve {\it 
admissible} if
it is surjective, and there is an enlargement $\pi \colon  \hat X \to X$
such that $f$ extends to a holomorphic map $\hat f\colon \hat X \to \bar C$
with connected fibers onto a smooth compactification $\bar C$ of $C$.
We will need the following generalization of the 
Beauville-Castelnuovo-De Franchis lemma \cite{B2, 2.1}.

\proclaim {Proposition 1.3} Let  $\bar L$ be  the Deligne extension of a unitary flat
line bundle $L= O_X\otimes \C_\psi$ on $X$ and $\theta \in H^0(\OmD {1} )$. 
Suppose that the sequence
$$ H^0(\bar X, \bar L) \lblarr {\theta} H^0(\bar X, \OmDL {1} ) \lblarr {\theta}
 H^0(\bar X, \OmDL {2} ) $$
is not exact. Then 
there exists an admissible holomorphic map $ f\colon  X \to C$  onto a
nonsimply connected smooth curve such that
 $\psi \in \rho f^*H^1(C, \C^*)$, where $\rho $ is torsion.
After enlarging $\bar X$ so that $f$ extends to $\bar f\colon  \bar X \to \bar C$,
we also have 
$ \theta \in \bar f^*H^0(\bar C , \Omega_{\bar C}^1(\bar C - C))$.
\endproclaim 
 
\pf The argument is very similar to the proof of \cite{B2, 2.1}.
By hypothesis, there is a nonzero form $\alpha \in H^0(\bar X, \OmDL {1} )$
such that $\theta \wedge \alpha = 0$; if $\bar L = O_{\bar X}$ then $\alpha$
can be chosen to not be a multiple of $\theta$. The relation $\theta \wedge
\alpha = 0$ implies that there exists a meromorphic section $\beta $ of 
$\bar L$ such that $ \alpha = \theta \otimes \beta$. By theorem 1.3 of 
chapter I, $\bar \del \alpha = 0$ and $d\theta = 0$, where $\bar \del $ is the 
unitary connection on $\bar L$. Hence $\theta \wedge \bar \del \beta = 0$.
This implies the existence of a meromorphic function $g$ on $\bar X$ such
that $\bar \del \beta = g \theta \otimes \beta$. Consequently
$$ dg \wedge \theta \otimes \beta = dg \wedge \theta \otimes \beta
+ g\theta \wedge \bar \del \beta = \bar \del ^2 \beta = 0.$$
Which implies $dg\wedge \theta = 0$.

 Suppose that $g$ were constant. Then $\alpha = \bar \del(\beta 
 /g)$. Since $\alpha$ is holomorphic, $\beta/g$ must holomorphic as
 well. This would imply that its derivative $\alpha$ must
 vanish (by 1.3 of chapter I); this contradicts the hypothesis. Therefore
$g$ is a nonconstant meromorphic map $\bar X -\,-> \PP^1$. Hironaka's
theorem implies that there is a nonsingular blow up $\pi\colon  \hat X \to \bar X$
such that $g\circ \pi$ extends to a holomorphic map $\hat g\colon  \hat X \to 
\PP^1$, and such that $\hat D = \pi^{-1}D$ is a divisor with normal 
crossings.
Let 
$ \hat X \lblarr {\hat f} \bar C \to \PP^1$
be the Stein factorization of  $\hat g$. The relation $dg \wedge \theta = 0$
implies that $\hat \theta =\pi^*\theta $ vanishes along the general fibers of $\hat f$.
Therefore $\hat \theta$ is the pullback of a meromorphic differential form from
$\bar C$ which must necessarily have logarithmic singularities.  
So there exists a divisor $E\subset \bar C$, with $\hat f^{-1} E \subseteq 
\hat D$, and
a form $\theta' \in H^0(\bar C, \Omega_{\bar C}^1(log E))$ such that 
$\hat \theta = \hat f^*\theta'$.
The curve $C = \bar C - E$ is not simply connected because it carries nonzero
logarithmic 1-form $\theta'$.
 If $C$ is compact then it has positive genus, and in
  this case the fibers of $\pi$ map trivially
to $C$, thus we may as well have taken $\hat X = \bar X$. In general, 
$\hat f(\pi^{-1}(X)) \subseteq C$ 
and we can assume equality holds after enlarging $E$.
So when $C$ is not compact, the fibers of the proper map 
$\pi|_{\pi^{-1}(X)}$ must map to points in $C$. Therefore 
$ \hat f |_{\pi^{-1}(X)}$ descends to a holomorphic map $f\colon X\to C$. Hence we can 
assume that $\hat X$ is an enlargement of $\bar X$, and we now replace $\bar X$
by $\hat X$ and $\theta$ by $\hat \theta$. 

Finally we have to show that $\psi\in \rho f^* H^1(C, \C^*)$ with $\rho$ a torsion
character. Let $F$ be the general fiber of $\hat f^{-1} (C) \to C$. Then
$im[H_1(F, \Z ) \to H_1(X, \Z)]$
has finite index in
$ ker[H_1(X, \Z) \to H_1(C, \Z)]$
by \cite{Se,1.3}. 
Therefore 
$ ker [H^1(X, \C^*) \to H^1(F, \C^*)]$
is a product of $f^*H^1(C, \C^*)$ and a torsion subgroup.
Thus it suffices to show that the restriction of $\psi$ to
$F$ is trivial. As $\alpha$ is nonzero, it will not vanish identically along
a general fiber $F$. Furthermore 
 $\theta = \hat f^*\theta'$ will not have any zeros along
$F$. Thus the relation $\alpha = \theta \otimes \beta$, implies
that $\beta$ is nonzero and holomorphic along $F$. Consequently
$$H^0(F\cap X, \C_\psi) \cong H^0(F, \bar L|_F) \not= 0,$$
which forces $\psi|_F$ to be trivial.

\proclaim {Proposition 1.4} If  $\xi \in \Sigma^1(X)$ is  nonunitary, then
there is an admissible map $f\colon X\to C$ onto a nonsimply connected curve
and a torsion character $\rho$ such that
$\xi \in \rho f^*H^1(C,\C^*)$.
\endproclaim 

\pf Let $N$ be the union of components about which $\xi$ has nontrivial monodromy.
After replacing $X$ by $\bar X - N$, we can assume that $\xi$ has nontrivial
local  monodromy about all components of $D$. 
As explained in  chapter II, 
$$H^1( X ,\C) = (H^{10}(X) \oplus H_\R^{11}(X)) \oplus  \sqrt{-1}(H_\R^{11}(X)
\oplus H^1(X, \R)).$$
The first summand is the subspace of $H^0(\bar X,\OmD {1} )$ consisting of
forms with purely imaginary residues. Therefore $\xi$ can be decomposed as
a product $\xi= \exp(\theta'+ \phi + r)$, with  $\theta' \in H^{10}+ H_\R^{11}$,
$\phi \in \sqrt{-1}H_\R^{11 }$ and $r \in \sqrt{-1}H^1(X,\R)$. 
After replacing $\phi + r$ with translate by an element of $H^1(X,\Z(1)) =
2\pi\sqrt{-1}H^1(X,\Z)$ we can assume that all the residues of $\phi$ are
nontrivial. Set 
$\theta = \theta' + \epsilon \phi$ and $\psi = \exp((1-\epsilon )\phi+r)$.
For suitable  real $\epsilon$, we can assume that $\theta$ is a nonzero element
of $H^0(\bar X,\OmD {1} )$ with nonintegral residues,  and  $\psi$ is a unitary
character with nontrivial monodromy about all components of $D$.

Let $L = O_X\otimes \C_\psi $ with flat unitary connection $\del$, and
let $(\bar L, \bar \del) $ denote the Deligne extension to $\bar X$. The monodromy
of the flat connection $\del -\theta$ is $\xi$ by lemma 1.1 of chapter I.
Therefore by proposition 1.2 of chapter I,
$$  {\HH}^1(\bar X, \OmL ; \bar \del + \theta) \cong H^1(X,\C_\xi) \not= 0.$$
There is a   spectral sequence converging to the left hand side with
$$ E_1^{pq} = H^q(\bar X, \OmDL {p} )$$
with differential $d_1$ induced by $\theta\wedge $ (because $\bar \del$ 
acts trivially on $E_1$). 
Therefore
either $E_2^{10}\not= 0$ or $E_2^{01} \not= 0$. The condition 
$E_2^{10}\not= 0$ is exactly the hypothesis of proposition 1.3, thus
after enlarging $\bar X$ 
there is an admissible map $f \colon  X \to C$ and a torsion character $\rho$
such that $\psi\in \rho f^*H^1(C,\C^*)$ and $\theta\in f^*H^0(\bar C, 
\Omega_{\bar C}^1(\log f(D)))$. This implies that $\xi \in \rho 
f^*H^1(C,\C^*)$.  

Now suppose that $E_2^{01} \not=0$. This means that
there is a nonzero class $\alpha \in H^1(\bar X, \bar L)$ such that $\theta \wedge 
\alpha \not= 0$.  
By assumption $(L,\del)$ has nontrivial monodromy about the components 
of $D$. This implies 
$$W_0(\Omega_{\bar X}^p(\log D)\otimes L^{\pm 1}) =
\Omega_{\bar X}^p(\log D)\otimes L^{\pm 1}. $$
By results of Timmerscheidt \cite{T1} (see also chapter IV, section 1),  
$H^q(\bar X, \OmDL {p} )$ is isomorphic to the space of harmonic $(p,q)$ forms
on $X$ satisfying $L^2$ growth conditions with respect to a suitable complete
finite volume K\"ahler metric on $ X$.
Therefore $\bar \alpha\in H^0(\bar X,\Omega_{\bar X}^1(\log D)\otimes {\bar L}^{-1})$ and 
$\theta \wedge\bar \alpha \in H^0(\bar X, \Omega_{\bar X}^1(\log D)\otimes {\bar L}^{-1}) $
 are $L^2$ forms. 
Let $\omega $ be the K\"ahler form, then the norm
$$ || \theta \wedge\bar \alpha ||^2 = 
\int_X \theta\wedge \bar \alpha\wedge \bar \theta\wedge \alpha \wedge 
\omega^{dimX-2} = 0$$
Thus $\theta \wedge \bar \alpha = 0$. Furthermore when $\bar L = O_{\bar X}$ then 
$ \alpha \not= c\bar\theta$ with $c\in \C^*$, since otherwise the norm 
$$||\theta||^2 = \int_{X} \theta\wedge\bar \theta \wedge \omega^{dimX-1} = 0.$$
Therefore proposition 1.3 again shows that after enlarging $\bar X$
there is an admissible map $f \colon  X \to C$ and a torsion character $\rho$
such that $\psi^{-1}\in \rho f^*H^1(C,\C^*)$ and $-\theta\in f^*H^0(\bar C, 
\Omega_{\bar C}^1(\log f(D)))$. Thus $\xi \in \rho^{-1}f^*H^1(C,\C^*)$.
\bigskip

Let us call two  admissible maps $f\colon  X \to C$ and $f' \colon  X \to C'$ equivalent
if there is an isomorphism $\sigma \colon  C\to C'$ such that $f' = f \circ \sigma$.

\proclaim {Lemma 1.5} The set $\calA (X)$ of equivalence classes of
 admissible maps $f\colon X\to C$, with
 $C$ nonsimply connected, is at most countable.
\endproclaim 

\pf Any such map $f\colon X \to C$ is determined by a surjective homomorphism
$Alb(X) \to Alb(C)$, and there are only countably many such maps.

\proclaim {Theorem 1.6} There exists a finite number of torsion characters 
$\rho_i \in H^1(X, \C)$,
unitary characters $\rho_j'$
and admissible holomorphic maps onto smooth curves $f_i\colon X\to C_i$ such
that
$$\Sigma^1(X) = \bigcup_i \rho_i f_i ^*H^1(C_i, \C^*) \,\cup \, \bigcup_j \rho_j$$
\endproclaim 

\pf  Let us write
$H(f,\rho) = \rho f^*H^1(f(X),\C^*).$ 
Any isolated point of $\Sigma^1(X)$ is unitary by proposition 1.4.
Now suppose that  $V$ is a positive dimensional irreducible component of $\Sigma^1(X)$.
By proposition 1.4, $V$ is contained in the union of 
$H(f,\rho)$ with $f\in \calA (X)$ and $\rho \in H^1(X,\C^*)_{tors}$. As an irreducible
variety cannot be a countable union of proper subvarieties, it follows that
$V $ is contained in one of the varieties $H(f,\rho)$. And we 
will show that they are  equal.
If $C= f(X)$ is $\C^*$,  then $dim\, H^1(C,\C^*) = 1$ 
so that $V =H(f,\rho)$. If $C$ is an elliptic curve, then 
 $V \subseteq H(f,\rho) \subseteq H^1(\bar X, \C^*)$
so $V$ will be the exponential of the translate of nonzero
sub Hodge structure of $H^1(C)$ by remark 2 above. Once again this
forces  $V =H(f,\rho)$.
In all other
cases the topological  Euler characteristic of $C$ is negative, therefore
$H^1(C,\C_{\xi}) \not= 0$ for all $\xi \in H^1(C,\C^*)$.
As the map 
$ H^1(C,\C_{\xi}) \to H^1(X, \C_{f^*\xi})$
is injective, $H(f,1) \subseteq \Sigma^1(X)$. Thus $V=H(f,\rho)$ when $\rho =1$.

 The remaining case when $C$ has negative Euler    characteristic and $\rho 
 \not= 1 $ is similar. We first ``untwist'' $\rho$.
As $\rho$ is torsion,  it lifts to a class $\tilde \rho\in H^1(X,\mu_N)$ 
 where $\mu_N \cong \Z/N\Z$ is the group of $N$th roots of unity for some 
 integer $N$. Let $p\colon Y\to X$ be a torser corresponding to $\tilde \rho$; it is
 an unramified cover with Galois group $\mu_N$, such that 
$p^*\rho = 1$.  From the exact sequence 
$$ 0 \to Hom(\mu_N,\C^*) \to Hom(\pi_1(X),\C^*)\to Hom(\pi_1(Y),\C^*) $$
we deduce that $\rho$ corresponds to a character of $\mu_N$. 
$p_*\C_Y$  decomposes, under the natural $\mu_N$ action,
 into a sum of one dimensional
local systems indexed by characters of $\mu_N$.
Therefore $\mu_N$ acts naturally on 
 $$H^1(Y,\C_{p^*\psi}) \cong H^1(X, p_*\C_Y \otimes \C_\psi )$$
  for any 
$\psi \in H^1(X,\C^*)$,  and $H^1(X, \C_{\rho\psi})$ is the $\rho$th
isotypic component of $H^1(Y,\C_{p^*\psi})$.
Let $\dim_\rho$ denote the dimension of the $\rho$th isotypic component
of a $\mu_N$-module.
To prove the theorem, it will suffice to  show that
$$dim_\rho H^1(Y, \C_{p^*f^*\xi}) \not= 0$$ 
for any $\xi \in H^1(C,\C^*)$. Let 
$ Y \lblarr {g} B \lblarr {\pi} C$
be a ``Stein factorization'' of $f\circ p$, i.e. a Stein factorization of a
fiberwise compactification. Then $\pi$ is a (possibly branched) $\mu_N$ covering
and $g$ is $\mu_N$ equivariant. Choose  finite triangulation $\Delta(C)$ of $C$, so 
that the branch points of $\pi$ are among the vertices. Let $\Delta^i(C)$ represent
the set of $i$-simplices.
 Lift $\Delta(C)$ to
a $\mu_N$ invariant triangulation $\Delta(B)$ of $B$. Let $\xi \in H^1(C,\C^*)$ and
let $S^\dt$ be the simplicial cochain complex associated to $\Delta(B)$ with
coefficients in  the local system $\C_{\pi^*\xi}$.
By additivity, the Euler characteristic
$$e_\rho = \sum_i\, (-1)^i dim_\rho  H^i(B,\C_{\pi^*\xi}) = \sum_i\, dim_\rho S^i.$$
As $\mu_N$ acts freely on the edges and faces, $S^i$is 
a free $\C[\mu_N]$-module of rank $\#\Delta^i(C)$, for $i>0$. $S^0$ is a direct sum 
$$ \oplus_{x\in \Delta_0(C)}\, \C[\mu_N/G_x]$$
where $G_x$ is the isotropy group of an element of $\pi^{-1}(x)$.
This implies that $dim_\rho(S^i) \le \#\Delta^i(C)$ with equality when
$i= 1,2$. Therefore $e_\rho < 0 $ which implies 
$dim_\rho H^1(B,\C_{\pi^*\xi}) \not= 0$. As $H^i(B,\C_{\pi^*\xi})$ is a $\mu_N$-submodule
of $H^1(Y, \C_{p^*f^*\xi})$, the theorem is proved.
\bigskip

 It is now
a straight forward exercise to generalize results concerning maps
of compact K\"ahler to curves (\cite{A}, \cite{B2},
\cite{C}, \cite{GL2}), to noncompact manifolds. 
We will say that a smooth curve is of {\it general type} if its topological
Euler characteristic is negative. (This is equivalent to the ampleness
of the log canonical bundle, so the terminology is consistent with
standard usage.) For such a curve, $H^1(C, \C_\rho) \not= 0$ for all
characters $\rho$. Given an integer $b$, let $N_b(X)$ be the
cardinality
of the set of equivalence classes of admissible maps to curves of
general
type with first betti number $b$. 

\proclaim {Proposition 1.7} If $f\colon  X \to C$ is an admissible map to a
curve of general type, then $f^*H^1(C, \C^*)$ is a component
of $\Sigma^1(X)$. Conversely, any positive dimensional component
of $\Sigma^1(X)$  containing $1$ is of this form.
\endproclaim 

\pf Given $f\colon  X \to C$  satisfying the above hypothesis, $f^*H^1(C, \C^*)$
is contained in a component of $\Sigma^1(X)$ of the form $g^*H^1(B, \C^*)$
for some admissible map $g\colon  X \to B$. This implies that there is a factorization
$ Alb(X) \to Alb(B) \to Alb(C),$
and consequently a factorization
$ X \to B \to C.$
As $f$ has connected fibers, $B \cong C$. 

Conversely,  any positive dimensional component
of $\Sigma^1(X)$ is given by   $f^*H^1(C, \C^*)$ with
$f\colon  X \to C$ admissible. We have to show that $C$ has general type. 
$C$ is of not of general type if and only if  it is an elliptic curve or
$\PP^1$ with at most 2 points removed. These examples will  be ruled out
by showing that $\Sigma^1(C)$ is infinite. It will
suffice to prove that 
$$H^1(X,\C_{f^*\rho}) \cong H^1(C,\C_\rho)$$
 for infinitely many $\rho$, or equivalently (by the Leray spectral
 sequence ) that $H^0(C,R^1f_*\C\otimes \C_\rho) = 0$ for infinitely
 many $\rho$. There is a Zariski open subset $j\colon U \hookrightarrow C$
 such that $f^{-1}U \to U$ is a locally trivial fibration. We claim
 that the canonical map $R^1f_*\C \to j_*j^*R^1f_*\C$ is injective.
 The injectivity can be checked at stalks around points of $C-U$.
Choose a small disk $\Delta$ centered  $p\in C-U$, and let $q\in \Delta-\{p\}$.
Then $(R^1f_*\C)_p \cong H^1(f^{-1}(\Delta),\C)$ and
$(j_*j^*R^1f_*\C)_p$ is isomorphic the monodromy invariant part
of $H^1(f^{-1}(q),\C)$. Consider the commutative diagram:
$$\matrix
                  &    & H^1(f^{-1}(\Delta      ))&\lblarr {r} 
                                                         & H^1(f^{-1}(q))&\cr
                  &    &  i\downarrow                 &   & \downarrow &\cr
H^1(\Delta -\{p\})& \to& H^1(f^{-1}(\Delta-\{p\}))&\to& H^1(f^{-1}(q))&\cr
\endmatrix
$$
Straightforward topological arguments show that bottom row is exact and
$i$ is injective. As $f$ is surjective, we can find a continuous cross
section $\sigma \colon  \Delta \to f^{-1}(\Delta)$. The map from $H^1(f^{-1}(\Delta))$
to $H^1(\Delta-\{p\})$ factors through $H^1(\Delta)$.
Thus the image of $H^1(f^{-1}(\Delta))$ 
in $ H^1(f^{-1}(\Delta-\{p\}))$ is contained in $ker[\sigma |_{\Delta-\{p\}}^*]$.
Therefore $r$
is injective and this implies the claim. Therefore
$$H^0(C,R^1f_*\C\otimes \C_\rho) = H^0(U,R^1f_*\C\otimes \C_\rho).$$
We can regard the local system $j^*R^1f_*\C$ as a
$\pi_1(U)$-module. Let $g_1,\ldots g_n$ be generators for $\pi_1(U)$. 
They act on the local system $j^*R^1f_*\C$ through linear transformations
$T_1,\ldots T_n$. Then
$$H^0(U, R^1f_*\C\otimes\C_\rho) \not= 0$$
only if $\rho(g_i)^{-1}$ are eigenvalues of $T_i$. Therefore the 
above group vanishes for almost all $\rho$.

\proclaim {Corollary 1.8} The number $N_b(X)$ is finite and depends
only on the fundamental group. More precisely, if $Y$ is a Zariski
open subset of a compact K\"ahler manifold with 
$\pi_1(Y) \cong \pi_1(X)$, then $N_b(Y) = N_b(X)$.
\endproclaim 

\pf The set 
$$\Sigma^1(X) = \Sigma^1(\pi_1(X))=\{\rho \in Hom(\pi_1(X), \C^*)|
H^1(\pi_1(X),\C_\rho) \not= 0\}$$
clearly depends only on $\pi_1(X)$, and $N_b(X)$ is just
the number of $b$ dimensional components containing $1$.

\proclaim {Corollary 1.9} $X$ possesses an admissible
map on a smooth curve of general type if and only if
$\pi_1(X)$ surjects homomorphically onto a nonabelian free group.
\endproclaim 

\pf The fundamental group of curve of general type is
either a nonabelian free group or is isomorphic to 
$$ < a_1,\ldots a_g,b_1,\ldots b_g | [a_1,b_1]\ldots[a_g,b_g]=1> $$
with $g \ge 2$. The latter group surjects onto the free
generated by $a_1,\ldots a_g$. This proves one direction.

Conversely, if $\pi_1(X)$ surjects onto a nonabelian free group
$F$ then $\Sigma^1(X) \supseteq \Sigma^1(F)$. Therefore $\Sigma^1(X)$
has an infinite component containing $1$.
\bigskip

Let $G' = [G,G]$ denote the commutator subgroup of a group $G$.
Recall that the first homology of $G$ with coefficients in
the trivial $G$-module $\C$ is $H_1(G,\C) \cong G/G'\otimes \C$ 

\proclaim {Corollary 1.10} If $X$ possesses an admissible map to curve
of general type then $dim H_1(\pi_1(X)', \C) = \infty$. If
 $dim H_1(\pi_1(X)', \C) = \infty$ then there is a finite unramified
abelian covering $Y \to X$, such that $Y$ possesses an admissible map
to a curve of general type.
\endproclaim 

\pf Let $M = H_1(\pi_1(X)',\C)$ and $\Lambda = H_1(X,\Z)
=\pi_1(X)/\pi_1(X)'$.
 $M$ has a natural $\Lambda$ action coming form the
group extension
$$1 \to \pi_1(X)' \to \pi_1(X) \to \Lambda \to 0.$$
In fact, $M$  is a finitely generated $\C[\Lambda]$ module;
thus it determines a coherent algebraic sheaf $\cal M$ on the
character variety $H^1(X,\C^*) = spec\C[\Lambda](\C)$ such
that $M = H^0({\cal M})$.
Using the Hochschild-Serre spectral sequence, one  can show that
$H^1(X,\C_\rho) \cong Hom_{\Lambda}(M, \C_\rho)$
when $\rho \not= 1$ (see \cite{A1 pp. 541-542} or \cite{B2, pp. 12-13}).
Therefore 
$$\Sigma^1(X) = supp {\cal M} \cup \{1\}$$
when $X$ has nonzero first betti number. Thus $\Sigma^1(X)$ is 
infinite if and only if $dim M =\infty$. $\Sigma^1(X)$ is infinite
if $X$ maps onto a curve of general type. 
On the other hand, by
theorem 1.6, if $\Sigma^1(X)$ is infinite then there is an
abelian cover of $Y\to X$ which has an infinite component containing
$1$. Therefore $Y$ possesses an admissible map onto a curve of general type.


\widestnumber\key{GL2}

\Refs

\ref \key A1
 \by D. Arapura
 \paper Hodge theory with local coefficients on
compact varieties
\jour Duke Math. J 
\vol 61 \yr 1990 \pages 531-543
\endref

\ref \key A
\bysame
\paper Higgs line bundles, Green-Lazarsfeld sets
and maps of K\"ahler manifolds to curves 
\jour Bull. AMS \yr 1992 \pages 310-314
\endref

\ref\key B1
\by A. Beauville
\paper Annulation du $H^1$ et syst\'emes
paracanoniques sur les surfaces
\jour J. reine angew. Math.
\vol 388 \yr 1988 \pages 149-157
\endref

\ref \key B2
\bysame
\paper Annulation du $H^1$ pour les fibr\'es en
en droites plats
\book Lect. Notes in Math. 1507
\publ Springer-Verlag \yr 1992 \pages 1-15
\endref

\ref\key C 
\by F. Catanese
\paper Moduli and classication of
irregular Kaehler manifolds (and algebraic varieties) with Albanese
general type fibrations
\jour Inv. Math. \vol 104 \yr 1991 \pages 263-289
\endref



\ref \key D1 
\by P. Deligne
\book Equation differentielles a point
singular regulier , Lect. Notes in Math. 163,
\publ Springer-Verlag \yr 1969
\endref

\ref\key D \bysame
\paper Theorie de Hodge II, III
\jour Publ. IHES 40,44
 \yr1972, 1974, \pages 5-57, 6-77
\endref

\ref\key G
\by R. Godement
\book Topologie alg\'ebrique et th\'eorie de faiseux 
\publ Hermman \yr1958
\endref

\ref\key GL1
\by M. Green, R. Lazarsfeld,
\paper Deformation theory, generic
vanishing theorems and some
conjectures of Enriques, Catanese and Beauville 
\jour Inv. Math. \vol 90
\yr 1987 \pages 389-407
\endref


\ref\key GL2 \bysame
\paper Higher obstructions to
deforming cohomology groups of line bundles
\jour J. A.M.S.
\vol 4 \yr1991 pages 87-103
\endref


\ref\key LM
\by A. Lubotzky, A. Magid
\paper Varieties of representations of
finitely generated groups
\jour Mem. AMS 336 \yr 1985
\endref

\ref\key MF
\by D. Mumford, J. Forgarty
\book Geometric invariant theory 
\publ Springer-Verlag \yr 1980
\endref

\ref\key N
\by V. Navaro Aznar
\paper Sur le theorie de Hodge-Deligne
\jour Inv. Math. \vol 90 \yr 1987 \pages 11-76
\endref



\ref\key Sa1
\by M. Saito
\paper Mixed Hodge modules
\jour Publ. Res. Inst. Math. Sci.
\vol 26 \yr 1990 \pages 221--333
\endref

\ref\key Sa2
\bysame
\paper Mixed Hodge modules and Admissible variations
 \jour CR Acad. \vol 309 \yr 1989 \pages  351-356
\endref

\ref\key Se
\by  F. Serrano
\paper Multiple fibers of a morphism
\jour  Comm. Math. Helv. \vol 65 \yr 1990 \pages 287-298
\endref


\ref\key S1
\by C. Simpson
\paper Higgs bundles and local systems
\jour Publ. IHES \vol 75 \yr 1992 \pages 5-95 
\endref


\ref\key S2
\bysame 
\paper Subspaces of moduli spaces of rank one local systems 
\jour Anal. ENS \vol 26 \yr 1993 \pages 361-401
\endref


\ref\key S3
\bysame 
\paper Moduli of representations of the fundamental group
of a projective variety II
\jour Publ. IHES \vol 80 \yr 1994\pages 5-79
\endref


\ref\key T1
\by K. Timmerscheidt
\paper Hodge theory for unitary local systems
\jour Inv. Math. \vol 86 \yr 1986 \pages 189-194
\endref


\ref\key T
\bysame
\paper  Mixed Hodge structures associated to unitary
 local systems
 \jour J. reine angew. Math. \vol 379 \yr 1987\pages 152-171
\endref


\ref\key Z
\by S. Zucker
\paper Hodge theory with degenerating coefficients
\jour Ann. Math. \vol 109 \yr 1979 \pages 415-476
\endref


\endRefs
\enddocument